\documentclass[onefignum,onetabnum]{siamart171218}   

\headers{Multilinear compressive sensing}{F. Malgouyres and J. Landsberg}

\title{Multilinear compressive sensing and an application to convolutional linear networks\thanks{Submitted to the editors 7/3/2018.}}

\author{Fran\c{c}ois Malgouyres
\thanks{Institut de Math\'ematiques de Toulouse ; UMR5219
Universit\'e de Toulouse ; CNRS UPS IMT; F-31062 Toulouse Cedex 9, France; (\email{Francois.Malgouyres@math.univ-toulouse.fr}) | and | Institut de Recherche Technologique Saint Exupery
}
\funding{This work has been supported by the DEEL program on Dependable and Explainable Learning (\url{www.deel.ai})}
\and Joseph Landsberg
\thanks{Department of Mathematics; Mailstop 3368;Texas A\& M University;
College Station, TX 77843-3368; (\email{jml@math.tamu.edu})}
\funding{Joseph Landsberg is supported by NSF DMS-1405348 and AF1814254}
}

\usepackage{hyperref}
\usepackage{amsmath}
\usepackage{amsfonts}
\usepackage{amsbsy}
\usepackage{amssymb}
\usepackage{relsize}
\usepackage{graphicx}
\usepackage{calc}
\usepackage{amstext}
\usepackage{bbm}
\usepackage{mathrsfs}
\usepackage{color}

\usepackage{ifthen}
\newboolean{long}   
\setboolean{long}{true}   

\usepackage{pict2e}
\linethickness{0.2mm}

\usepackage{graphicx}

\usepackage{xcolor}

\usepackage{algorithm}
\usepackage{algorithmic}

\def\endproof{\hfill $\Box$\newline\newline}  

\DeclareMathOperator {\supp} {Supp}
\newcommand {\SUPP}[1] {\supp\left( #1 \right)}
\newcommand {\one} {\mathbbm{1}}
\newcommand {\RR} {\mathbb R}
\newcommand {\CC} {\mathbb C}

\newcommand {\NN} {\mathbb N}
\newcommand {\NNN}[1] {[#1]}

\newcommand {\calA} {\mathcal A}

\newcommand {\calS} {\mathcal S}
\newcommand {\calT} {\mathcal T}

\newcommand {\calF} {\mathcal{F}}
\newcommand {\cald} {\delta}

\newcommand {\Mod}[1] {{\mathcal M}^{#1}}
\newcommand {\TS} {\RR^{S^K}}   
\newcommand {\hS} {\RR^{S\times K}}    
\newcommand {\hSdiag} {\RR^{S\times K}_{\mbox{\tiny diag}}}    
\newcommand {\iS} {\NNN{S}^{K}}  
\newcommand {\ibf} { \mathbf{i}}      
\newcommand {\Ibf} { \mathbf{I}}      
\newcommand {\jbf} { \mathbf{j}}      
\newcommand {\ebf} { \mathbf{e}}      
\newcommand {\hbf} { \mathbf{h} }      
\newcommand {\gbf} { \mathbf{g}}      
\newcommand {\pbf} { \mathbf{p}}      

\newcommand {\PS}[2] {\langle #1 , #2 \rangle}     
\DeclareMathOperator {\rk} {rk}
\newcommand {\RK}[1] {\rk\left( #1 \right)}   
\DeclareMathOperator {\myspan} {Span}
\newcommand {\SPAN}[1] { \myspan\left( #1 \right) }
\DeclareMathOperator {\myker} {Ker}
\newcommand {\KER}[1] { \myker\left( #1 \right) }

\newcommand {\DIM}[1] {\dim\left( #1 \right)}   
\DeclareMathOperator {\argmin} {argmin}
\DeclareMathOperator {\argmax} {argmax}
\DeclareMathOperator {\diag} {diag}
\newcommand {\DIAG}[1] { \diag \left( #1 \right) }   
\newcommand {\CL}[1] { \langle #1 \rangle }   

\newcommand {\tree}{\mathcal G}
\newcommand {\nodes} {\mathcal N}
\newcommand {\edges} {\mathcal E}
\newcommand {\leaves} {\mathcal F}

\newcommand {\PP} {{\mathcal P}}
\newcommand {\multiconv}[2] { {\mathcal T}^{#1}(#2)\, }

\newcommand {\codeset} {\mathbb R^{N  |\leaves|}}
\newcommand {\RP} {\mathbb R^N}

\newcommand {\class}[1] {\{ #1 \}}
\newcommand {\NSP} {deep-NSP}
\newcommand {\NSPlong} {deep-Null Space Property}

\newsiamthm{defi}{Definition}
\newsiamthm{prop}{Proposition}
\newsiamthm{thm}{Theorem}
\newsiamthm{cor}{Corollary}
\newsiamremark{com}{Comment}
\newsiamthm{TargStat}{Informal theorem}

\def\PROOFUN{

\subsection{Proof of  \cref{metric-prop}}\label{metric-prop-proof}

Notice that,  the sets $\CL{\hbf}\cap\hSdiag$ and $\CL{\gbf} \cap\hSdiag$ are finite and therefore the infimum in the definition of $d$ is reached. We also have whatever $\hbf$, $\gbf \in \hS_*$
\begin{equation}\label{lhbhlveajvgt}
d_p(\CL{\hbf}, \CL{\gbf}) =  \inf_{\hbf'\in\CL{\hbf}\cap\hSdiag}
\Bigl(\inf_{\gbf'\in\CL{\gbf} \cap\hSdiag  } \|\hbf'-\gbf'\|_p\Bigr).
\end{equation}

Moreover, whatever $\hbf \in \hS_*$ and $\hbf'$ and $\hbf'' \in \CL{\hbf} \cap\hSdiag$ there exist $(s_k)_{k\in\NNN{K}} \in\{-1,1\}^K$ such that $\prod_{k\in\NNN{K}} s_k =1$ and
\[ \hbf'_k = s_k \hbf''_k\qquad, \forall k\in\NNN{K}.
\]
Using the above two properties, we can check that
$$\inf_{\gbf'\in\CL{\gbf} \cap\hSdiag  } \|\hbf'-\gbf'\|_p = \inf_{\gbf'\in\CL{\gbf} \cap\hSdiag  } \|\hbf''-\gbf'\|_p$$
As a consequence, the outer infimum in \cref{lhbhlveajvgt} is irrelevant and we have
\[d_p(\CL{\hbf}, \CL{\gbf}) =  \inf_{\gbf'\in\CL{\gbf} \cap\hSdiag  } \|\hbf'-\gbf'\|_p\qquad,\forall \hbf, \, \gbf \in \hS_* \mbox{ and }\hbf' \in \CL{\hbf} \cap\hSdiag.\]

Using this last property, we easily check that $d_p$ is a metric on $\hS_* /\sim$.
\endproof
}

\def\PROOFDEUX{

\subsection{Proof of  \cref{rk1-Ident-thm}}\label{rk1-Ident-thm-proof}

Notice first that when $K=1$  the inequality is a straightforward consequence of the usual inequalities between $l^p$ norms. We therefore assume from now on that $K\geq 2$.

All along the proof, we consider $\hbf$ and $\gbf\in\hS_*$ and assume that $\|P(\hbf)\|_\infty \geq \|P(\gbf)\|_\infty $. We also assume that $\|P(\gbf)-P(\hbf)\|_\infty \leq \frac{1}{2}\|P(\hbf)\|_\infty $. We first prove the inequality when $p=q=+\infty$.

In order to do so, we consider
\[\ibf\in \argmax_{\jbf \in \iS}~|P(\hbf)_\jbf|
\]
and assume, without lost of generality (otherwise, we can multiply one vector of $\hbf$ and $\gbf$ by $-1$ to get this property and multiply back once the inequality have been established), that $P(\hbf)_\ibf\geq 0$. We therefore have $P(\hbf)_\ibf = \|P(\hbf)\|_\infty$. Notice also that we have, under the above hypotheses,
\begin{equation}\label{qribnqervn}
\|P(\gbf)\|_\infty \geq P(\gbf)_\ibf\geq  P(\hbf)_\ibf - \|P(\gbf)-P(\hbf)\|_\infty \geq \frac{1}{2}\|P(\hbf)\|_\infty >0.
\end{equation}

Moreover, we consider the operator $E_\ibf$ that extracts the $K$ signals of size $S$ that are obtained when freezing,  at the index $\ibf$ in a tensor $T$, all coordinates but one. Formally, we denote
\begin{eqnarray*}
E_\ibf:\TS & \longrightarrow & \hS \\
T & \longmapsto & E_\ibf(T)
\end{eqnarray*}
where for all $k\in\NNN{K}$ and all $j\in\NNN{S}$
\[E_\ibf(T)_{k,j} = T_{\ibf_1,\ldots,\ibf_{k-1},j,\ibf_{k+1},\ldots,\ibf_K}.
\]
We consider
\[\hbf' = (P(\hbf)_\ibf)^{-1+\frac{1}{K}} ~E_\ibf( P(\hbf) )\qquad\mbox{and}\qquad\gbf' = (P(\gbf)_\ibf)^{-1+\frac{1}{K}}~ E_\ibf( P(\gbf) ).
\]
We have for all $\jbf\in\iS$
\begin{eqnarray*}
P(\hbf')_\jbf & = & (P(\hbf)_\ibf)^{-K+1} ~ P\left( E_\ibf( P(\hbf) )\right)_\jbf, \\
 & = & (P(\hbf)_\ibf)^{-K+1} ~ \prod_{k=1}^K P(\hbf)_{\ibf_1,\ldots,\ibf_{k-1},\jbf_k,\ibf_{k+1},\ldots,\ibf_K}\\
 & = & (P(\hbf)_\ibf)^{-K+1} ~ \prod_{k=1}^K \hbf_{1,\ibf_1} \ldots \hbf_{k-1,\ibf_{k-1}}
 \hbf_{k,\jbf_k} \hbf_{k+1,\ibf_{k+1}}\ldots \hbf_{K,\ibf_{K}}\\
& = &   \prod_{k=1}^K \hbf_{k,\jbf_k} = P(\hbf)_\jbf.
\end{eqnarray*}
We therefore have $ P(\hbf') =  P(\hbf)$. This can be written $\hbf'\in\CL{\hbf}$. Similarly, we have $ \gbf' \in \CL{\gbf}$.

Also, because of the definition of $\ibf$ and $\hbf'$, we are guaranteed that, whatever $k\in\NNN{K}$,
\begin{eqnarray*}
\|\hbf'_k\|_\infty & = & (P(\hbf)_\ibf)^{-1+\frac{1}{K}} \|E_\ibf( P(\hbf) )\|_\infty \\
& = &\|P(\hbf)\|_\infty^{-1+\frac{1}{K}} \|P(\hbf)\|_\infty = \|P(\hbf)\|_\infty^{\frac{1}{K}}
\end{eqnarray*}
The latter being independent of $k$, we have $\hbf'\in\hSdiag$. Unfortunately, unless for instance $\ibf\in\argmax_{\jbf\in \iS} | P(\gbf)_\jbf| $, it might occur that $\gbf'\not\in\hSdiag$. However, if we consider
\[\gbf'' \in \argmin_{{\bf f}\in \CL{\gbf} \cap \hSdiag} \|{\bf f} - \gbf'\|_\infty,
\]
we have since $\hbf'\in \CL{\hbf}\cap\hSdiag $ and $\gbf''\in \CL{\gbf}\cap\hSdiag $
\[d_\infty(\CL{\hbf}, \CL{\gbf}) \leq \|\hbf' - \gbf''\|_\infty 
\]
and therefore
\begin{equation}\label{junbjtrsb}
 d_\infty(\CL{\hbf}, \CL{\gbf}) \leq \|\hbf' - \gbf'\|_\infty + \|\gbf' - \gbf''\|_\infty .
\end{equation}
In the sequel we will successively calculate upper bounds of $\|\hbf' - \gbf'\|_\infty$ and $\|\gbf' - \gbf''\|_\infty$ in order to find an upper bound of $d_{\infty}(\CL{\hbf}, \CL{\gbf})$.

\vspace*{0.5cm}
\noindent{\bf Upper bound of $\|\hbf' - \gbf'\|_\infty$:}

We have
\begin{eqnarray*}
\| \hbf'-\gbf' \|_\infty & = & \|  (P(\hbf)_\ibf)^{-1+\frac{1}{K}} ~E_\ibf( P(\hbf) ) - (P(\gbf)_\ibf)^{-1+\frac{1}{K}}~ E_\ibf( P(\gbf) ) \|_\infty \\
& \leq & \|  (P(\hbf)_\ibf)^{-1+\frac{1}{K}} \left( E_\ibf( P(\hbf) )-E_\ibf( P(\gbf) ) \right)\|_\infty \\
&& \qquad\qquad+ \|\left( (P(\hbf)_\ibf)^{-1+\frac{1}{K}} - (P(\gbf)_\ibf)^{-1+\frac{1}{K}} \right) E_\ibf( P(\gbf) ) \|_\infty \\
& \leq &\|P(\hbf)\|_\infty^{-1+\frac{1}{K}} \|E_\ibf( P(\hbf) )-E_\ibf( P(\gbf))\|_\infty + \|P(\gbf)\|_\infty |(P(\hbf)_\ibf)^{-1+\frac{1}{K}} - (P(\gbf)_\ibf)^{-1+\frac{1}{K}} |\\
& \leq &\|P(\hbf)\|_\infty^{-1+\frac{1}{K}} \| P(\hbf)- P(\gbf)\|_\infty + \|P(\hbf)\|_\infty |(P(\hbf)_\ibf)^{-1+\frac{1}{K}} - (P(\gbf)_\ibf)^{-1+\frac{1}{K}} |
\end{eqnarray*}
But we also have using the mean value theorem and \cref{qribnqervn}
\begin{eqnarray*}
|(P(\hbf)_\ibf)^{-1+\frac{1}{K}} - (P(\gbf)_\ibf)^{-1+\frac{1}{K}} | &  \leq& \left(1-\frac{1}{K}\right) P(\gbf)_\ibf^{-2+\frac{1}{K}} |P(\hbf)_\ibf-P(\gbf)_\ibf| \\
&  \leq& \left(1-\frac{1}{K}\right) \left(\frac{1}{2} \|P(\hbf)\|_\infty\right)^{-2+\frac{1}{K}} \|P(\hbf)-P(\gbf)\|_\infty\\
&  \leq&4~\|P(\hbf)\|_\infty^{-2+\frac{1}{K}}\|P(\hbf)-P(\gbf)\|_\infty
\end{eqnarray*}
We therefore finally obtain that
\begin{equation}\label{ernontbtg}
\|\hbf' - \gbf'\|_\infty  \leq  5 \|P(\hbf)\|_\infty^{-1+\frac{1}{K}} \| P(\hbf)- P(\gbf)\|_\infty.
\end{equation}

\vspace*{0.5cm}
\noindent{\bf Upper bound of $\|\gbf' - \gbf''\|_\infty$:}

First, since $\gbf''\in \CL{\gbf} = \CL{\gbf'} $, we know that there exists $(\lambda_k)_{k\in\NNN{K}} \in×\RR^K$ such that
\begin{equation}\label{oerjv}
\prod_{k=1}^K \lambda_k = 1
\end{equation} and
\[\gbf''_k = \lambda_k \gbf'_k \qquad \mbox{, for all } k\in\NNN{K}.
\]
Furthermore, we have for all $k\in\NNN{K}$
\begin{equation}\label{eriipqzrniov}
\|\gbf'_k - \gbf''_k\|_\infty = |1-\lambda_k |~ \|\gbf'_k\|_\infty.
\end{equation}
Also, if there is $k'$ such that $\lambda_{k'} <0$, since \cref{oerjv} holds, there necessarily exist another $k''$ such that $\lambda_{k''} <0$. If we replace $\gbf''_{k'}$ by $-\gbf''_{k'}$ and replace $\gbf''_{k''}$ by $-\gbf''_{k''}$ we remain in $\CL{\gbf}\cap\hSdiag$ and can only make $\|\gbf' - \gbf''\|_\infty$ decrease. Repeating this process until all the $\lambda_k$'s are non-negative, we can assume without loss of generality that
\[ \lambda_k\geq0\qquad\mbox{, whatever } k\in\NNN{K}.
\]

This being said, we establish two other simple facts that motivate the structure of the proof. First, in order to find an upper bound for \cref{eriipqzrniov},  we easily establish (using \cref{qribnqervn}) that
\begin{eqnarray*}
\|\gbf'_k\|_\infty & = &  (P(\gbf)_\ibf)^{-1+\frac{1}{K}}~ \|E_\ibf( P(\gbf) ) \|_\infty  \\
& \leq  &  (\frac{1}{2} \|P(\hbf)\|_\infty)^{-1+\frac{1}{K}}~\|P(\hbf) \|_\infty 
\end{eqnarray*}
and therefore
\begin{equation}\label{iqeurhiohqe}
\|\gbf'_k\|_\infty \leq   2  \|P(\hbf)\|_\infty^{\frac{1}{K}} .
\end{equation}
Second, the value $\lambda_k$ appearing in \cref{eriipqzrniov}, can be bounded by using bounds on  $\|\gbf'_k\|_\infty$ and the identity
\begin{equation}\label{erqiqzrrf}
\|\gbf''_k\|_\infty = \|P(\gbf)\|_\infty^{\frac{1}{K}} = \lambda_k ~\|\gbf'_k\|_\infty.
\end{equation}
Qualitatively, the latter identity indeed guarantees that, as $\|P(\gbf) -P(\hbf)\|_\infty$ goes to $0$, $\lambda_k$ goes to $1$. Let us now establish this quantitatively.

Recalling   that
\[\gbf'=   (P(\gbf)_\ibf)^{-1+\frac{1}{K}}~ E_\ibf( P(\gbf) ),
\]
and using \cref{qribnqervn} again, we obtain
\[\|\gbf'_k\|_\infty \leq \left( \|P(\hbf) \|_\infty -\frac{1}{2} \|P(\hbf) - P(\gbf) \|_\infty  \right)^{-1+\frac{1}{K}} ~ \|P(\gbf) \|_\infty.
\]
We also have (again, using \cref{qribnqervn})
\begin{eqnarray*}
\|\gbf'_k\|_\infty & \geq & (P(\gbf)_\ibf)^{-1+\frac{1}{K}} |P(\gbf)_\ibf| \\
& = & (P(\gbf)_\ibf)^{\frac{1}{K}} \\
& \geq & \left(\|P(\hbf) \|_\infty -\frac{1}{2} \|P(\hbf) - P(\gbf) \|_\infty  \right)^{\frac{1}{K}}.
\end{eqnarray*}

Plugging the upper bound of $\|\gbf'_k\|_\infty$ in \cref{erqiqzrrf}, using successively \cref{qribnqervn}, the mean value theorem and the hypothesis on the size of $P(\hbf) - P(\gbf)$  gives:
\begin{eqnarray*}
\lambda_k - 1 & = &\frac{ \|P(\gbf)\|_\infty^{\frac{1}{K}} }{\|\gbf'_k\|_\infty} ~-1\\
& \geq & \|P(\gbf)\|_\infty^{-1+\frac{1}{K}}  \left( \|P(\hbf) \|_\infty -\frac{1}{2} \|P(\hbf) - P(\gbf) \|_\infty  \right)^{1-\frac{1}{K}} ~ -1 \\
& \geq &\left( 1 -\frac{\|P(\hbf) - P(\gbf) \|_\infty}{2\|P(\hbf) \|_\infty}   \right)^{1-\frac{1}{K}} ~ -1 \\
& \geq & - (1-\frac{1}{K}) \left( 1 -\frac{\|P(\hbf) - P(\gbf) \|_\infty}{2\|P(\hbf) \|_\infty}   \right)^{-\frac{1}{K}} ~ \frac{\|P(\hbf) - P(\gbf) \|_\infty}{2\|P(\hbf) \|_\infty} \\
& \geq & - \left(1 - \frac{1}{4}\right)^{-\frac{1}{K}} ~ \frac{\|P(\hbf) - P(\gbf) \|_\infty}{2\|P(\hbf) \|_\infty} \\
& \geq & - ~ \frac{\|P(\hbf) - P(\gbf) \|_\infty}{\|P(\hbf) \|_\infty}.
\end{eqnarray*}
Similarly, plugging the lower bound of $\|\gbf'_k\|_\infty$ in \cref{erqiqzrrf}, we obtain using successively \cref{qribnqervn}, the mean value theorem and the hypothesis on the size of $P(\hbf) - P(\gbf)$:

\begin{eqnarray*}
\lambda_k - 1 & \leq &  \|P(\gbf)\|_\infty^{\frac{1}{K}} \left(\|P(\hbf) \|_\infty -\frac{1}{2} \|P(\hbf) - P(\gbf) \|_\infty  \right)^{-\frac{1}{K}} ~ -1 \\
& \leq & \left( 1 -\frac{\|P(\hbf) - P(\gbf) \|_\infty}{2\|P(\hbf) \|_\infty}   \right)^{-\frac{1}{K}} ~ -1 \\
& \leq & \frac{1}{K} \left( 1 -\frac{\|P(\hbf) - P(\gbf) \|_\infty}{2\|P(\hbf) \|_\infty}   \right)^{-1-\frac{1}{K}} ~ \frac{\|P(\hbf) - P(\gbf) \|_\infty}{2\|P(\hbf) \|_\infty} \\
& \leq & \frac{1}{K} \left( 1 - \frac{1}{4}\right)^{-1-\frac{1}{K}} ~\frac{\|P(\hbf) - P(\gbf) \|_\infty}{2\|P(\hbf) \|_\infty} \\
&\leq & \frac{4^2}{2K3^2}~ \frac{\|P(\hbf) - P(\gbf) \|_\infty}{\|P(\hbf) \|_\infty} \\
& \leq & \frac{\|P(\hbf) - P(\gbf) \|_\infty}{\|P(\hbf) \|_\infty}.
\end{eqnarray*}

Finally, we get
\begin{equation}\label{rqgfniobg}
| \lambda_k - 1 | \leq  \frac{\|P(\hbf) - P(\gbf) \|_\infty}{\|P(\hbf) \|_\infty}.
\end{equation}

By combining \cref{eriipqzrniov}, \cref{iqeurhiohqe} and \cref{rqgfniobg},
we obtain
\[\|\gbf'_k - \gbf''_k\|_\infty \leq 2 ~  \|P(\hbf)\|_\infty^{-1+\frac{1}{K}} ~ \|P(\hbf) - P(\gbf) \|_\infty .
\]

Combining the latter inequality with \cref{junbjtrsb} and \cref{ernontbtg} provides
\[d_\infty(\CL{\hbf}, \CL{\gbf}) \leq 7 \|P(\hbf)\|_\infty^{-1+\frac{1}{K}} \| P(\hbf)- P(\gbf)\|_\infty,
\]
and concludes the proof when $p=q=+\infty$.

In order to establish the property when $1\leq p \leq +\infty$ and $1\leq q \leq +\infty$, we simply use the fact that
\[d_p(\CL{\hbf}, \CL{\gbf})\leq (KS)^{\frac{1}{p}} ~ d_\infty(\CL{\hbf}, \CL{\gbf})
\]
and
\[\| P(\hbf)- P(\gbf)\|_\infty \leq \| P(\hbf)- P(\gbf)\|_q.
\]

\endproof
}
\def\PROOFTROIS{

\subsection{Proof of \cref{rk1-Ident-sharp}}\label{rk1-Ident-sharp-proof}

In the example, we consider $\hbf$ and $\gbf$ such that for all $k\in\NNN{K}$ and all $i\in\NNN{S}$
\[\hbf_{k,i} = \left\{\begin{array}{ll}
1 & \mbox{ if } i=0, \\
0 & \mbox{ otherwise,}
\end{array}\right. \qquad \mbox{ and } \qquad
\gbf_{k,i} =\left\{\begin{array}{ll}
\left(\frac{1}{2}\right)^{\frac{1}{K}} & \mbox{ if } i=0, \\
\epsilon_q & \mbox{ otherwise,}
\end{array}\right.
\]
where $\epsilon_{+\infty} = \left(\frac{1}{2}\right)^{\frac{1}{K}}$ and $\epsilon_q =\min\left(  \left(\frac{ 1-\left(\frac{1}{2}\right)^{ \frac{q}{K} } }{ S-1 }\right)^{\frac{1}{q}} , \left(\frac{1}{2}\right)^{\frac{1}{K}} \right)$, if $q<+\infty$. We immediately obtain
\[\|P(\hbf) \|_\infty= 1,\qquad\|P(\gbf) \|_\infty= \frac{1}{2}\qquad\mbox{and}\qquad\|P(\hbf)- P(\gbf) \|_\infty= \frac{1}{2}.
\]
We also have,
\[ d_p(\CL{\hbf},\CL{\gbf})^p = \| \hbf -\gbf \|_p^p \geq K(S-1)~\epsilon_q^{p}\geq \frac{KS}{2}~\epsilon_q^{p}.
\]
Decomposing the sum necessary to the calculation of the $l^q$ norm of a tensor according to number of index different from $0$ (which corresponds to $l$ in the sum below), we obtain
\begin{eqnarray*}
\|P(\hbf)- P(\gbf) \|_q^q & = & \sum_{l=0}^K \binom{l}{K} (S-1)^l \epsilon_q^{lq} \left(\frac{1}{2}\right)^{ \frac{(K-l)q}{K} }, \\
 & = &\left( \left(\frac{1}{2}\right)^{ \frac{q}{K}}   + (S-1) \epsilon_q^q   \right)^K \leq 1.
\end{eqnarray*}

We then easily obtain that
\[
7 \|P(\hbf)\|_\infty^{-1+\frac{1}{K}} (KS)^{\frac{1}{p}} \| P(\hbf)- P(\gbf)\|_q   \leq  7 (KS)^{\frac{1}{p}}, 
\]
and therefore
\begin{equation} \label{qemibnq}
 7 \|P(\hbf)\|_\infty^{-1+\frac{1}{K}} (KS)^{\frac{1}{p}} \| P(\hbf)- P(\gbf)\|_q \leq  7  ~ \frac{ d_p(\CL{\hbf},\CL{\gbf}) }{\epsilon_q} ~ 2^{\frac{1}{p}}. 
\end{equation}
We first calculate a lower bound of $\epsilon_q$ when $\epsilon_q= \left(\frac{ 1-\left(\frac{1}{2}\right)^{ \frac{q}{K} } }{ S-1 }\right)^{\frac{1}{q}}$ (which, in particular, rules out $q=+\infty$). Using the mean value theorem, we obtain
\[1-\left(\frac{1}{2}\right)^{ \frac{q}{K}  } \geq \min_{t\in[\frac{1}{2},1]} \left(\frac{q}{K} t^{\frac{q}{K} - 1}\right)  (1-\frac{1}{2}).
\]
Distinguishing, whether $q\leq K$ or not, we find after a short calculation that, since $q\geq 1$,
\[1-\left(\frac{1}{2}\right)^{ \frac{q}{K}  } \geq \min \left( \frac{1}{2K}, \frac{1}{K} \left(\frac{1}{2}\right)^\frac{q}{K} \right) = \frac{1}{K} \min\left(\left(\frac{1}{2}\right)^\frac{1}{q},\left(\frac{1}{2}\right)^\frac{1}{K}\right)^q \geq \frac{1}{K2^q}.
\]
We then deduce
\[\epsilon_q \geq\frac{1}{2 \left(K S \right)^{\frac{1}{q}} }.
\]
Of course, when $\epsilon_q=\left(\frac{1}{2}\right)^{\frac{1}{K}}$ (which includes $q=+\infty$), we  immediately obtain
\[\epsilon_q \geq \frac{1}{2}.
\]
Using this lower bound in \cref{qemibnq} leads to the bounds stated in the proposition.
\endproof

}

\def\PROOFQUATRE{

\subsection{Proof of \cref{PLip-thm}}\label{PLip-thm-proof}

Before starting the proof, we define for any $k\in\{0,\ldots,K\}$
\[P_k(\hbf,\gbf)_\ibf = \gbf_{1,\ibf_1}\ldots\gbf_{k,\ibf_k}\hbf_{k+1,\ibf_{k+1}}\ldots\hbf_{K,\ibf_K}\qquad \mbox{, for all } \hbf, \gbf\in\hS \mbox{ and all }\ibf\in\iS.
\]

We consider $\gbf$ and $\hbf\in\hS$. Let us first assume that $\|\gbf\|_\infty \leq \|\hbf\|_\infty = 1 $. We have for any $\ibf\in\iS$, using this hypothesis and standard inequalities between $l^p$ norms, when $q<+\infty$
\begin{eqnarray*}
|P(\gbf)_\ibf - P(\hbf)_\ibf |^q & = & \left| \sum_{k=0}^{K-1} \left( P_{k+1}(\hbf,\gbf)_\ibf - P_k(\hbf,\gbf)_\ibf\right) \right|^q \\
& \leq & K^{q-1}  \sum_{k=0}^{K-1}\left|P_{k+1}(\hbf,\gbf)_\ibf - P_k(\hbf,\gbf)_\ibf\right|^q \\
& \leq & K^{q-1}  \sum_{k=0}^{K-1} \left|\gbf_{k+1,\ibf_{k+1}} -\hbf_{k+1,\ibf_{k+1}} \right|^q
\end{eqnarray*}
The same calculation when $q=+\infty$ leads to
\[|P(\gbf)_\ibf - P(\hbf)_\ibf | \leq K \max_{k=1..K}\left|\gbf_{k,\ibf_{k}} -\hbf_{k,\ibf_{k}}\right|.
\]

Therefore, we have when $q<+\infty$
\begin{eqnarray*}
\|P(\hbf) - P(\gbf) \|^q_q & = & \sum_{\ibf \in \iS}  |P(\hbf)_\ibf - P(\gbf)_\ibf|^q \\
& \leq & K^{q-1}  \sum_{k=1}^{K} \sum_{\ibf \in \iS}  \left|\gbf_{k,\ibf_{k}} -\hbf_{k,\ibf_{k}} \right|^q \\
& = & K^{q-1}  \sum_{k=1}^{K}  S^{K-1}  \|\gbf_k - \hbf_k\|_q^q\\
& =& K^{q-1} S^{K-1} \|\gbf - \hbf\|_q^q
\end{eqnarray*}
and therefore
\[\|P(\hbf) - P(\gbf) \|_q \leq K^{1-\frac{1}{q}} S^{\frac{K-1}{q}} \|\gbf - \hbf\|_q.
\]
Again, a similar calculus for $q=+\infty$ leads to
\[\|P(\hbf) - P(\gbf) \|_{+\infty} \leq K \|\gbf - \hbf\|_{+\infty}.
\]
Remember that the two last inequalities hold for $\gbf$ and $\hbf\in\hS$ such that $\|\gbf\|_\infty \leq \|\hbf\|_\infty = 1 $.

Let us now consider any $\gbf'$ and $\hbf'\in\hS$ and any $\gbf\in\hSdiag\cap \CL{\gbf'}$ and $\hbf\in\hSdiag\cap \CL{\hbf'}$. We denote $\delta = \max(\|\gbf\|_{+\infty}, \|\hbf\|_{+\infty})$. Notice first that $\|\gbf\|_{+\infty} = \|P(\gbf') \|_{+\infty} ^{\frac{1}{K}}$ and $\|\hbf\|_{+\infty} = \|P(\hbf') \|_{+\infty} ^{\frac{1}{K}}$. Therefore
\begin{equation}\label{etbvontb}
\delta = \max(\|P(\gbf') \|_{+\infty} , \|P(\hbf') \|_{+\infty} )^\frac{1}{K}.
\end{equation}

We can apply the above inequality to $\frac{\hbf}{\delta}$ and $\frac{\gbf}{\delta}$ (we might need to switch $\hbf$ and $\gbf$ but it does not change the final inequality) and obtain when $q<+\infty$
\[\|P(\frac{\hbf}{\delta}) - P(\frac{\gbf}{\delta}) \|_q \leq K^{1-\frac{1}{q}} S^{\frac{K-1}{q}} \|\frac{\gbf}{\delta} - \frac{\hbf}{\delta}\|_q.
\]
This leads to
\[\|P(\hbf) - P(\gbf) \|_q \leq K^{1-\frac{1}{q}} S^{\frac{K-1}{q}} \delta^{K-1} \|\gbf - \hbf\|_q.
\]
Similarly,  when $q=+\infty$, we obtain
\[\|P(\hbf) - P(\gbf) \|_{+\infty} \leq K \delta^{K-1} \|\gbf - \hbf\|_{+\infty}.
\]
The fact that these two last inequalities hold for any $\gbf\in\hSdiag \cap \CL{\gbf'}$ and any $\hbf\in\hSdiag\cap \CL{\hbf'}$, together with \cref{etbvontb}, leads to the statement provided in  \cref{PLip-thm}.
\endproof
}

\def\PROOFCINQ{

\subsection{Proof of  \cref{rkA-prop}}\label{rkA-prop-proof}

 The span of the Segre variety $P(\hS)$ is the full ambient space $\TS$, so there exists sets of
$R\leq S^K$ points on it that are linearly independent.
The set of $R$-tuples of points on  $P(\hS)$
that fail to be linearly independent  is a proper subvariety
of the variety of sets of $R$-tuples of points on $P(\hS)$
because being a linearly independent set of points  is an open condition and there
exists sets of points that are linearly independent. Therefore
$R\leq S^K$ independent and randomly chosen points according to a continuous distribution on $P(\hS)$ will be linearly independent.

The intersection $P(\hS) \cap \KER{\calA}$ is
a proper subvariety of $P(\hS)$, so with probability one,
$R\leq S^K$ independent randomly chosen points according to a continuous distribution will not intersect it and be linearly
independent. This is indeed the intersection of two non-empty open conditions. Therefore, all spans of subsets of the points will
intersect $\KER{\calA}$ transversely (in particular, the span of fewer than $\RK{\calA}$ points
will not intersect it). Thus there image under $\calA$ will have dimension as large
as possible. The same argument works if $R>S^K$.
\endproof

}
\def\PROOFSIX{
\subsection{Proof of  \cref{lifting-prop}}\label{lifting-Prop-proof}

The proof relies on the fact that for any $T^* \in \argmin_{T\in\TS}  \|\calA T - X \|^2$,  we have
\[\calA^t (\calA T^*-X) = 0,
\]
where $\calA^t:  \RR^{n\times m}\rightarrow  \TS$ is the adjoint  linear map.  This implies that
 for any $$T^* \in \argmin_{T\in\TS}  \|\calA T - X \|^2,$$ any $L\in\NN$ and any $\hbf\in\Mod{L}$
\begin{eqnarray*}
\| \calA P(\hbf) - X \|^2 & = &  \| \calA(P(\hbf)-T^* ) + (\calA T^*- X) \|^2,\nonumber \\
	 & = &  \| \calA(P(\hbf)-T^* ) \|^2 +  \|\calA T^*- X \|^2  + 2 \PS{\calA(P(\hbf)-T^* )}{\calA T^*- X}, \nonumber\\
	 & = &  \| \calA(P(\hbf)-T^* ) \|^2 +  \|\calA T^*- X \|^2.
\end{eqnarray*}
In words, $\| \calA P(\hbf) - X \|^2$ and $\| \calA(P(\hbf)-T^* ) \|^2$ only differ by an additive constant. Moreover, since the value of the objective function  $\|\calA T^*- X \|^2$ is independent of  the particular minimizer $T^*$ we are considering, this additive constant is independent of $T^*$. As a consequence,
a minimizer of $\| \calA P(\hbf) - X \|^2$ also minimizes $\| \calA(P(\hbf)-T^* ) \|^2$ and vice versa.
\endproof

}

\def\PROOFSEPT{

\subsection{Proof of \cref{carac-thm}}\label{carac-thm-proof}

 Write  $\overline T = P(\overline \hbf)$ and let $L^*$ and $\hbf^*$ be a minimizer of \cref{model}.  \Cref{lifting-prop} and the fact that $\overline T$ minimizes \cref{Pb1} implies that $(L^*,\hbf^*)\in \argmin_{L\in\NN, \hbf\in\Mod{L}} \| \calA(P(\hbf)- \overline T) \|^2$. As a consequence,
\[\| \calA(P(\hbf^*)-\overline{T} ) \|^2 = 0
\]
and
\[P(\hbf^*) \in \overline{T} +\KER{\calA},
\]
proving the first implication.

Conversely, let  $L^*\in\NN$ and $\hbf^*\in\Mod{L^*}$
be  such that $P(\hbf^*)\in \overline{T} + \KER{\calA}$, then
\[\|\calA(P(\hbf^*)-\overline{T}) \|^2= 0 = \min_{L\in\NN, \hbf\in\Mod{L}}  \|\calA (P(\hbf) -\overline{T}) \|^2.
\]
As a consequence, $(L^*,\hbf^*) \in\argmin_{L\in\NN, \hbf\in\Mod{L}}  \|\calA (P(\hbf) -\overline{T}) \|^2$ and, using \cref{lifting-prop}, $\hbf^*$ is a minimizer of \cref{model}.
\endproof

}

\def\PROOFHUIT{

\subsection{Proof of  \cref{ident-thm}}\label{ident-thm-proof}
\begin{itemize}
\item  Proof of the first statement of  \cref{ident-thm}:

We first assume that $\CL{\overline\hbf}$ is identifiable. We consider $L^*$ and $h^*$ such that there is $L^*$ such that $P(\hbf^*) \in (P(\overline\hbf)+\KER{\calA})\cap P(\Mod{L^*})$.  We know from \cref{carac-thm} that $\hbf^*\in \argmin_{L\in\NN,\hbf\in\Mod{L}} \| \calA P(\hbf)-X \|^2$. Using that $\CL{\overline\hbf}$ is identifiable,  $\CL{\hbf^*}=\CL{\overline\hbf}$ and, from Standard Fact  \cref{rk1-Ident-cor} (at the beginning of \cref{tensor-sec}), we get $P(\hbf^*) = P(\overline\hbf)$. Finally, we can conclude, that if $\CL{\overline\hbf}$ is identifiable we have $(P(\overline\hbf)+\KER{\calA})\cap P(\Mod{}) \subset \{P(\overline\hbf)\}$.

Let us assume now that for all $L\in\NN$, $(P(\overline\hbf)+\KER{\calA})\cap P(\Mod{L}) \subset \{P(\overline\hbf)\}$ and consider $$(L^*,\hbf^*)\in \argmin_{L\in\NN, \hbf\in\Mod{}} \| \calA P(\hbf)-X \|^2.$$ Using  \cref{carac-thm}, we know that $P(\hbf^*) \in (P(\overline\hbf)+\KER{\calA})\cap P(\Mod{L^*})$. Using the hypothesis, we have $P(\hbf^*)=P(\overline\hbf)$ and using Standard Fact \cref{rk1-Ident-cor}, we finally conclude that $\CL{\hbf^*}=\CL{\overline\hbf}$. This completes the proof of the first statement.

\item Proof of the second statement of  \cref{ident-thm}:

 Assume that there is $L$ and $L'\in\NN$ such that $\KER{\calA}\cap \bigl(P(\Mod{L})-P(\Mod{L'}) \bigr) \not\subset \{0\}$ then there exist $\hbf\in\Mod{L}$ and $\overline \hbf\in\Mod{L'}$ such that $P(\hbf)\neq P(\overline{\hbf})$ and $P(\hbf)-P(\overline{\hbf})\in\KER{\calA}$. Using the first statement of the proposition, we obtain that $\overline{\hbf}$  is not identifiable. As a conclusion, $\Mod{}$ is not identifiable.

Conversely, assume that there exists $L'$ and some non-identifiable  $\overline{\hbf}\in\Mod{L'}$. Using the first statement of the proposition, we know that there exists $L\in\NN$ and $\hbf\in\Mod{L}$ such that $P(\hbf)\neq P(\overline{\hbf})$ and $P(\hbf)- P(\overline{\hbf}) \in\KER{\calA}$. Therefore $\KER{\calA}\cap\bigl(P(\Mod{})-P(\Mod{}) \bigr) \not\subset \{0\}$.

\end{itemize}
\endproof

}

\def\PROOFNINE{

\subsection{Proof of  \cref{suffic-ident-thm}}\label{suffic-ident-thm-proof}

We first make the \lq\lq equality holds  generically\rq\rq\ statement precise in
our context. Fix any variety $X$ and assume $Y$ is a linear space, say
of dimension $y$. Let $G(y,\CC^N)$ denote the Grassmannian of $y$-planes
through the origin in $\CC^N$. The Grassmannian is both a smooth manifold
and an algebraic variety. We can interpret \lq\lq equality holds  generically\rq\rq\  in this
context as saying for a Zariski open subset of $G(y,\CC^N)$, equality will
hold. In our situation, if we fix $\rk(\calA)$ and allow $\ker(\calA)$ to vary
as a point in the Grassmannian, with probability one, it will intersect
$J(P(\Mod{L}), P(\Mod{L'}))$ only in the origin, and this assertion
is also true over $\RR$ because complex numbers are only needed
to assure existence of intersections, not non-existence.
\endproof

}

\def\PROOFDIX{

\subsection{Proof of \cref{suf-stble-recovery}}\label{suf-stble-recovery-proof}

We have
\begin{eqnarray*}
\|\calA (P(\hbf^*) - P(\overline\hbf) )\| & \leq & \|\calA P(\hbf^*) - X \| + \|\calA P(\overline\hbf) -X\| \\
& \leq &  \delta + \eta
\end{eqnarray*}

Geometrically, this means that $P(\hbf^*)$ belongs to a cylinder centered at $P(\overline\hbf)$ whose direction is $\KER{\calA}$ and whose section is defined using the operator $\calA$. If we further decompose (the decomposition is unique)
\[P(\hbf^*) - P(\overline\hbf) = T + T',
\]
where $T'\in\KER{\calA}$ and $T$ is orthogonal to $\KER{\calA}$, we have
\begin{equation}\label{sigmamin}
\|\calA (P(\hbf^*) - P(\overline\hbf) )\| = \|\calA T\| \geq \sigma_{min} \|T\|,
\end{equation}
where $\sigma_{min}$ is the smallest non-zero singular value of $\calA$. We finally obtain
\[
 \|P(\hbf^*) - P(\overline\hbf) -T' \| = \|T\| \leq \frac{\delta+\eta}{\sigma_{min}}.
\]
The term on the left-hand side corresponds to the distance between a point in $P(\Mod{L^*}) - P(\Mod{\overline{L}})$ (namely $P(\hbf^*) - P(\overline\hbf)$) and a point in $\KER{\calA}$ (namely $T'$).

Since $\KER{\calA}$ satisfies the \NSP~with constants $(\gamma,\rho)$, when  $ \delta+\eta \leq \rho$, we obtain the first inequality of the theorem
\[ \|P(\hbf^*) - P(\overline\hbf)\| \leq \gamma ~\frac{\delta+\eta}{\sigma_{min}}.
\]
When $\overline{\hbf} \in\hS_*$, for $\frac{\gamma}{ \sigma_{min}} ~(\delta+\eta) \leq \frac{1}{2}  ~ \max\left(\|P(\overline\hbf)\|_\infty,\|P(\hbf^*)\|_\infty\right)$, we can apply \cref{rk1-Ident-thm} and obtain \eqref{tkl,stb}.
\endproof

}

\def\PROOFONZE{

\subsection{Proof of  \cref{nec-stble-recovery}}\label{nec-stble-recovery-proof}

Let $\overline{L}$ and $\overline{L}'\in\NN$ and $\overline{\hbf}\in\Mod{\overline{L}}$ and $\overline{\hbf}'\in\Mod{\overline{L}'}$ be such that $\|\calA\left(P(\overline{\hbf})- P(\overline{\hbf}')\right)\|\leq \delta$. We also consider throughout the proof $T'\in\KER{\calA}$. We assume that $\|P(\overline{\hbf})\|_\infty \leq \|P(\overline{\hbf}')\|_\infty$. If it is not the case, we simply switch  $\overline{\hbf}$ and $\overline{\hbf}'$ in the definition of $X$ and $e$ below. We denote
\[X=\calA P(\overline{\hbf})\qquad\mbox{ and } \qquad e=\calA P(\overline{\hbf})-\calA P(\overline{\hbf}').
\]
We have $X=\calA P(\overline{\hbf}')+e$ and $\|e\|\leq \delta$. Therefore, the hypothesis of the theorem (applied with $\hbf^* = \overline\hbf$ and $L^*=\overline L$) guarantees that
\[d_2(\CL{\overline\hbf} , \CL{\overline\hbf'})\leq C ~ \|P(\overline\hbf')\|_\infty^{\frac{1}{K}-1} \|e\|.
\]

Using the fact that $e=\calA P(\overline{\hbf})-\calA P(\overline{\hbf}')$ and $T'\in\KER{\calA}$ we obtain
\[\|e\| = \|\calA ( P(\overline{\hbf}) -P(\overline{\hbf}') - T')  \| \leq \sigma_{max}~ \|P(\overline{\hbf})- P(\overline{\hbf}') - T'  \|.
\]
where $\sigma_{max}$ is the spectral radius of $\calA$. Therefore
\[d_2(\CL{\overline\hbf} , \CL{\overline\hbf'})\leq C \|P(\overline\hbf')\|_\infty^{\frac{1}{K}-1} ~\sigma_{max}~ \|P(\overline{\hbf})- P(\overline{\hbf}') - T'  \|.
\]

Finally, using  \cref{PLip-thm} and the fact that $\|P(\overline{\hbf})\|_\infty \leq \|P(\overline{\hbf}')\|_\infty$, we obtain
\begin{eqnarray*}
\|P(\overline\hbf')-P(\overline\hbf)\| & \leq & S^{\frac{K-1}{2}} K^{1-\frac{1}{2}}\|P(\overline\hbf')\|_{\infty}^{1-\frac{1}{K}} d_2(\CL{\overline\hbf'}, \CL{\overline\hbf}) \\
 & \leq & C S^{\frac{K-1}{2}} \sqrt{K} ~\sigma_{max}~\|P(\overline{\hbf})- P(\overline{\hbf}') - T'  \| \\
 & = & \gamma \|P(\overline{\hbf})- P(\overline{\hbf}') - T'  \|
\end{eqnarray*}
for $\gamma = C S^{\frac{K-1}{2}} \sqrt{K}~\sigma_{max}~$.

Summarizing, we conclude that under the hypothesis of the theorem: For any $T\in P(\Mod{\overline{L}}) -  P(\Mod{\overline{L}'}) $ such that $\|\calA T\|\leq \delta$ we have for any $T'\in\KER{\calA}$
\[ \|T\|\leq \gamma \| T - T' \|.
\]
\endproof

}

\def\PROOFDOUZE{

\subsection{Proof of \cref{nec-ident-network}}\label{nec-ident-network-proof}

Throughout the proof, we define, for any $\ibf\in\iS$, $\hbf^\ibf\in\hS$ by 
\begin{equation}\label{hibf}
\hbf^\ibf_{k,j} = \left\{\begin{array}{ll}
1 & \mbox{ if } j =\ibf_k \\
0 & \mbox{otherwise}
\end{array}\right. \qquad\mbox{, for all }k\in\NNN{K}\mbox{ and }  j\in\NNN{S}.
\end{equation}
This notation shall not be confused with $\hbf^\pbf$, with $\pbf\in\PP$.
\begin{itemize}
\item Let us first prove the first statement:
We can easily check that $(P(\hbf^\ibf))_{\ibf\not\in\Ibf}$ forms a basis of $\{T\in\TS \mid \forall \ibf\in\Ibf, T_\ibf = 0\}$. We can also easily check using \cref{multiconv} that, for any $\ibf\not\in\Ibf$,
\[\calA P(\hbf^\ibf) =  M_1(\hbf^\ibf_1)\ldots M_K(\hbf^\ibf_K) = 0.
\]
Therefore,  $\{T\in\TS \mid \forall \ibf\in\Ibf, T_\ibf = 0\} \subset \KER{\calA}$.

Conversely, for any  $\ibf\in\Ibf$, we can deduce from \cref{multiconv} and the hypotheses of the proposition that all the entries of $\calA P(\hbf^\ibf)$ are in $\{0,1\}$. We denote $D_\ibf = \{(i,j)\in\NNN{N}\times\NNN{N|\leaves|}  \mid  \calA P(\hbf^\ibf)_{i,j} = 1\}$. Using (again) the hypothesis of the proposition and \cref{multiconv}, we can prove that, for any distinct $\ibf$ and $\jbf\in\Ibf$, we have $D_\ibf \cap D_\jbf = \emptyset$. This easily leads to the \cref{premier_item} of the first statement. We also deduce that
\[\RK{\calA} \geq |\Ibf| = S^K - \dim(\{T\in\TS \mid \forall \ibf\in\Ibf, T_\ibf = 0\}).
\]
Finally, we deduce that $\dim(\KER{\calA}) \leq \dim(\{T\in\TS \mid \forall \ibf\in\Ibf, T_\ibf = 0\})$ and therefore 
\[\KER{\calA} = \{T\in\TS \mid \forall \ibf\in\Ibf, T_\ibf = 0\}.
\]

\item Let us now prove the second statement:
Using the hypothesis of the second statement and \cref{multiconv}, we know that there is $f\in\leaves$ and $n\in\NNN{N}$ such that
\[\sum_{p\in\PP(f)} \multiconv{p}{\one}_n \geq 2.
\]
As a consequence, there is $\ibf$ and $\jbf\in\iS$ with $\ibf\neq \jbf$ and
\[\multiconv{\pbf_\ibf}{\hbf^\ibf}_n = \multiconv{\pbf_\jbf}{\hbf^\jbf}_n = 1.
\]
Therefore,
$\calA P(\hbf^\ibf) = \calA P(\hbf^\jbf)
$
and the network is not identifiable.
\end{itemize}
\endproof

}

\def\PROOFTREISE{

\subsection{Proof of  \cref{network-cor}}\label{network-cor-proof}

The fact that, under the hypotheses of the proposition, $\KER{\calA} =\{0\}$ is a direct consequence of  \cref{nec-ident-network}. The \NSP~property and the value of $\gamma$ also follow from the definition of the \NSP. 

To calculate $\sigma_{min}$, let us consider $T\in\TS$ and express it under the form $T=\sum_{\ibf \in\Ibf} T_\ibf P(\hbf^\ibf)$, where $\hbf^\ibf$ is defined \cref{hibf}.
Let us also remind that, applying  \cref{nec-ident-network}, the supports of $\calA P(\hbf^\ibf)$ and $\calA P(\hbf^\jbf)$ are disjoint, when $\ibf\neq\jbf$. Let us finally add that, since $\calA P(\hbf^\jbf)$ is the matrix of a convolution with a Dirac mass, its support is of size $N$. We finally have
\begin{eqnarray*}
\|\calA T \|^2 & = & \|\sum_{\ibf \in\Ibf} T_\ibf \calA P(\hbf^\ibf)  \|^2, \\
& = & N \sum_{\ibf \in\Ibf} T_\ibf^2 = N\|T\|^2,
\end{eqnarray*}
from which we deduce the value of $\sigma_{min}$.
\endproof

}

\def\PROOFQUATORZE{

\subsection{Proof of  \cref{stabl-rec-network-thm}}\label{stabl-rec-network-thm-proof}

Considering  \cref{nec-ident-network}, we only need to prove that the condition is sufficient to guarantee the parameter stability.

Let us consider a path $\pbf\in\PP$, using \cref{multiconv}, since all the entries of $M_1(\one)\ldots M_K(\one)$ belong to $\{0,1\}$, all the entries of $M_1(\one^\pbf)\ldots M_K(\one^\pbf)$ belong to $\{0,1\}$. Therefore, we can apply \cref{network-cor} and  \cref{suf-stble-recovery} to the restriction of the convolutional linear network to $\pbf$ and obtain
\[d_p ( \CL{(\hbf^*)^\pbf} , \CL{\overline\hbf^\pbf} ) \leq \frac{7 (KS')^{\frac{1}{p}}}{\sqrt{N}}     \min\left( \|P(\overline\hbf^\pbf)\|_\infty^{\frac{1}{K}-1}, \|P((\hbf^*)^\pbf)\|_\infty^{\frac{1}{K}-1} \right) (\delta^\pbf+\eta^\pbf),
\] 
where $\delta^\pbf$ and $\eta^\pbf$ are the restrictions of the errors on $\SUPP{\pbf}$.

We therefore have
\[d_p ( \CL{(\hbf^*)^\pbf} , \CL{\overline\hbf^\pbf} ) \leq \frac{7 (KS')^{\frac{1}{p}}}{\sqrt{N}} \varepsilon^{1-K} (\delta^\pbf+\eta^\pbf),
\]
and finally
\begin{eqnarray*}
\cald_p (\class{\hbf^*}, \class{\overline\hbf}) & \leq & \frac{7 (KS')^{\frac{1}{p}} \varepsilon^{1-K}}{\sqrt{N}} \left(\sum_{\pbf\in \PP} (\delta^\pbf+\eta^\pbf)^p \right)^{\frac{1}{p}}, \\
& \leq & \frac{7 (KS')^{\frac{1}{p}} \varepsilon^{1-K}}{\sqrt{N}}(\delta+\eta).\\
\end{eqnarray*}
\endproof

}

\begin{document}

\maketitle

\begin{abstract}
  We study a deep linear network endowed with the following structure: a matrix $X$ is obtained by multiplying $K$ matrices (called factors and corresponding to the action of the layers). The action of each layer (i.e. factor) is obtained by applying a fixed linear operator to a vector of parameters satisfying a constraint. The number of layers is not limited. Assuming that $X$ is given and factors have been estimated, the error between the product of the estimated factors and $X$ (i.e. the reconstruction error) is either the statistical or the empirical risk. 
  
We provide necessary and sufficient conditions on the network topology under which a stability property holds. The stability property requires that the error on the parameters defining the near-optimal factors scales linearly with the reconstruction error (i.e. the risk). Therefore, under these conditions on the network topology, any successful learning task leads to stably defined features that can be interpreted.

In order to do so, we first evaluate how the Segre embedding and its inverse distort distances. Then, we show that  any deep structured linear network can be cast as a generic multilinear problem that uses the Segre embedding. This is the {\em tensorial lifting}. Using the tensorial lifting, we  provide a necessary and sufficient conditions for the identifiability of the factors up to a scale rearrangement. We finally provide  necessary and sufficient condition called \NSPlong~(because of the analogy with the usual Null Space Property in the compressed sensing framework) which guarantees that  the stability property holds. 

We illustrate the theory with a practical example where the deep structured linear network is a convolutional linear network. We obtain a condition on the scattering of the supports which is strong but not empty. A simple test on the network topology can be implemented to test if the condition holds.
\end{abstract}



\section{Introduction}
\subsection{The aim of the paper}

Deep learning has led to many practical breakthroughs and has led to significant improvements and state of the art performances in many fields such as computer vision, natural language processing, signal processing, robotics {\it etc}. The range of applications grows at a strong pace. Despite these empirical successes, the theory supporting deep learning is still far from satisfactory. For instance, sharp and accurate answers to the most natural questions on: the efficiency of optimization algorithms when applied to the objective function minimized in deep learning (\cite{livni2014computational,haeffele2015global,choromanska2015loss,choromanska2015open,safran2016quality,baldi1989neural,baldi1995learning,kawaguchi2016deep,2018arXiv180206384V}); and the expressiveness of the networks (\cite{barron1993universal,andoni2014learning,eldan2016power,cohen2016expressive,khrulkov2017expressive,cohen2016convolutional,montufar2014number,telgarsky2016benefits}); and guarantees on the statistical risk (\cite{janzamin2015beating,goel2017eigenvalue,xie2016diversity,schmidt2017nonparametric}) for learned neural networks are still missing. This makes it difficult to optimize and configure neural networks. Moreover, the absence of answers to these questions prevents certification that systems built with deep learning algorithms are robust. 

The reasons explaining the outcome of a neural network are often difficult to highlight (\cite{bach2015pixel,montavon2017explaining,hendricks2016generating}). Even worse, despite the settings described in \cite{Arora,arora2014provable,brutzkus2017globally,li2017convergence,sedghi2014provable,zhong2017recovery}, the instability of the parameters optimizing the deep learning objective does not allow the interpretation of the features defined by these parameters. This last problem is the one we investigate in this work.

Our goal in this paper is to evaluate how far the architectures used in applications are from architectures for which we can guarantee that the parameters returned by the algorithm, and therefore the features defined using these parameters, are stably defined. To do so, we consider two families of networks and establish necessary and sufficient conditions on their topology guaranteeing that the features learned by the algorithm are stably defined.

More precisely,  we establish statements of the following form for two families of deep networks. Below, the action of the network parameterized by $\hbf$ is denoted $f_\hbf$.

\begin{TargStat}{\bf Stability guarantee}\label{unformal-stat}

We assume a known parameterized family of functions $f_\hbf$ and a metric\footnote{The metric takes into account inter-layer rescaling.} $d$ between parameter pairs. We establish  a necessary and sufficient condition on the family $f_\hbf$ guaranteeing that:

There exists a constant $C>0$ such that for any input/output pairs $I$, $X$ and any pair of parameters $\hbf ^*$, $\overline{\hbf}$ for which
\[\delta = \|X- f_{\hbf^*}(I) \|,
\]
and
\[\eta = \|X- f_{\overline{\hbf}}(I) \|,
\]
are sufficiently small,  we have
\begin{equation}\label{erug}
d(\overline h,h^*) \leq C (\delta + \eta).
\end{equation}
\end{TargStat}

Considering a regression problem, the values $\delta$ and $\eta$ can be interpreted as the statistical or the empirical risk for the parameters $\overline h$ and $h^*$. The inequality \cref{erug} therefore guarantees that the set made of the parameters leading to a small risk has a small diameter. The features defined using such parameters are therefore stably defined. This seems to be the minimal condition allowing the interpretation of the features. The condition on the family of functions $f_\hbf$ is typically a condition on the topology of the network.

In \cref{unformal-stat}, $\overline h$ and $h^*$ might have different roles. For instance, if we know that the input/output pairs have been generated using a particular $\overline h$, possibly up to some error as modeled by $\delta$, then \cref{erug} guarantees that $h^*$ is close to $\overline h$ and provides a way to control the statistical risk.

The existing stability guarantees \cite{Arora,arora2014provable,brutzkus2017globally,li2017convergence,sedghi2014provable,zhong2017recovery,MalgouyresLandsbergITW}
consider this setting and describe both a network topology and an algorithm whose output $h^*$ is guaranteed to be close to $\overline h$. In this study, we do not make any assumption on the construction of $\overline h$ and $h^*$ and our objective is more modest. With regard to their objective, giving a necessary and sufficient condition of stability plays the same role as a complexity theory statement saying that a particular configuration is NP-hard. It rules out some network topologies.

Notice that, when $I$ and $X$ are such that it is possible to have $\delta=\eta=0$, the above stability guarantee
implies that the minimizer of the network objective function is unique. \Cref{necess-ident-thm} and Theorem  \cref{suffic-ident-thm} will show that this uniqueness condition is strongly related to the level of over-specification of the network. The simplified and intuitive statement is that optimal solutions of overspecified networks are not unique and are unstable. This explains the instability observed in applications. The theorems analyze this property in detail. This might be viewed as a negative result since overspecification is currently the main hypothesis of statements guaranteeing the success of the neural network optimization.

\subsection{The considered deep networks}
\subsubsection{Overview}
We consider two kinds of deep networks: A general family of deep structured linear networks\footnote{We call this family {\em deep structured linear networks} because the family is endowed  with tools to impose structures. We analyze the impact of the structure on the stability property. However, these tools might be used to define the usual deep linear network.} in Sections \cref{ident-sec} and \cref{stable-sec};  and a family of convolutional linear networks in \cref{conv-tree}. The formal statements for the deep structured linear networks are in Theorems \cref{suf-stble-recovery} and \cref{nec-stble-recovery}. The statements for convolutional linear networks are in \cref{stabl-rec-network-thm}. Below, we describe the deep structured linear networks.

\subsubsection{Deep structured linear networks}\label{struct-linear-network-sec}

The term deep linear network usually corresponds to fully-connected feed-forward networks, without bias, in which the activation function is the identity. In the general results described in this paper we consider deep linear networks and provide two means to enforce some structure to the network. As we  describe below, the structures can be used to include: feed-forward linear networks;  convolutional linear networks (as is done in \cref{conv-tree}); the action of a ReLU activation function; sparse networks; non-negative networks; and combinations of the above. The family also includes most matrix factorization problems.

We model a deep structured linear network as a product of matrices called factors. The factors  depend linearly on parameters in $\RR^S$, for $S\in\NN$.

More precisely, consider an arbitrary depth parameter $K\geq 1$. The number of layers is $K+1$ and the layers are enumerated in such a way that the layer receiving the input is $K+1$ and the layer returning the output is $1$. We consider sizes $m_1 \hdots   m_{K+1} \in \NN$,  write $m_1=m$, $m_{K+1}=n$. We consider,  for $k=1\hdots K$, the linear map
\begin{equation}\label{mkdef}
\begin{array}{rcl}
M_k: \RR^S &\longrightarrow &\RR^{m_k\times m_{k+1}} \\
h &\longmapsto & M_k(h) 
\end{array}
\end{equation}
Given some parameters, $h_1, \hdots, h_K\in\RR^S$, the action of the deep structured linear network is the product
\[M_1(h_1) \cdots M_K(h_K)
\]
The factor $M_K(h_K)$ might involve the inputs of the samples by considering $M_K(h_K) = M'_K(h_K)I$ for: a linear map $M'_K$; and for a matrix $I$ whose columns contain the inputs. Given outputs $X\in\RR^{m\times n}$, the optimization of the parameters $h_1, \hdots, h_K$ defining the network aims at getting 
\[M_1(h_1) \cdots M_K(h_K) \simeq X.
\]

To model feed-forward linear networks, the mappings $M_k$, $k=1\hdots K-1$ (and $M_K'$) construct the matrix by placing the entry of $h_k$ corresponding to an edge in the network in the corresponding entry in  $M_k(h_k)$.

For convolutional layers, $M_k$ and $M'_K$ concatenate convolution matrices\footnote{Depending on the situation: Toeplitz, block-Toeplitz, circulant or block-circulant matrices. The matrices often involve downsampling.} defined by a portion of the entries in $\hbf_k$. Each convolution matrix is at the location corresponding to a prescribed edge.

The main argument for studying deep structured linear networks is due to their strong connection to non-linear networks that uses the rectified linear unit (ReLU)\footnote{ReLU is the most common activation function.} activation function. We explain it in detail.  The action of the ReLU activation function at the layer $k$ treats every entry independently of the other entries and multiplies it by either $1$ (the entry is kept) or $0$ (the entry is canceled). More precisely, denoting $\hbf = (h_k)_{k=1..K}$, the action of the ReLU activation function on the layer $k$ is to apply the map
$A_k : \RR^{m_k\times n}\longmapsto \RR^{m_k\times n}$  (where $m_k\times n$ is the size data in the layer $k$) such that:
\[(A_kM)_{i,j} = a_k(\hbf)_{i,j} M_{i,j},\qquad \mbox{ for } (i,j) \in\{1,\hdots, m_k\}\times\{1,\hdots, n\}
\]
where $a_k(\hbf)\in\{0,1\}^{m_k\times n}$ is defined by
\[a_k(\hbf)_{i,j} = \left\{\begin{array}{ll}
1 & \mbox{ if } \Bigl(M_{k+1}(h_{k+1}) A_{k+1}M_{k+2}(h_{k+2}) \cdots A_{K-1}M_{K}(h_{K})\Bigr)_{i,j} \geq 0 \\
0 & \mbox{ otherwise.}
\end{array}\right.
\]
The function
\begin{eqnarray*}
a_k: \hS & \longrightarrow & \{0,1\}^{m_k\times n} \\
\hbf & \longmapsto & a_k(\hbf)
\end{eqnarray*}
is piecewise constant because $\{0,1\}^{m_k\times n}$ is finite. (This has already been used in \cite{safran2016quality}.) As a consequence, the parameter space $\hS$ is partitioned into subsets such that on every subset $a_k$ is constant, for all $k=1..K$. Therefore, on every subset the action of the non-linear network coincides with the action of a deep structured  linear network that groups at every layer $A_k$ and $M_{k+1}$. Further, the landscape of the objective function of the non-linear neural network that uses ReLU coincides, on every part of the partition, with the landscape of a deep structured linear network. This is a strong argument in favor of the study of deep structured linear networks.

Notice that deep structured linear networks have also been obtained in \cite{choromanska2015loss,choromanska2015open,kawaguchi2016deep} by modeling  the action of the activation function  as random, independent of the input and when considering the expectation of the network action. However, these assumptions are not satisfied by the deep networks used in applications (see \cite{choromanska2015open}) and it is not clear that this link can be exploited to obtain theoretical guarantees for realistic deep networks.

In addition to the structure induced by the operators $M_k$, we also consider
structure imposed  on the vectors $h$. We assume that we know a collection of models $\Mod{}=(\Mod{L})_{L\in\NN}$ with the property that for every $L$, $\Mod{L} \subset \hS$
is a given subset. We will assume that the parameters $\hbf\in\hS$ defining the factors are such that there exists $L\in\NN$ such that $\hbf\in\Mod{L}$. For instance, the constraint $\hbf\in\Mod{L}$ might be used to impose sparsity, grouped sparsity or co-sparsity. One might also use the constraint $\hbf\in\Mod{L}$ to impose non-negativity, orthogonality, equality, compactness, {\it etc}. Generally speaking $\Mod{}$ is used to impose some prior or some form of regularity or to compress the parameter space and obtain better bounds \cite{arora2018stronger}. The models might also be used to alleviate ambiguities. For instance, if the operators $M_k$ and $M_{k+1}$ allow permutations (i.e. there exists $(h_k,h_{k+1})\neq (g_k,g_{k+1})$ and a permutation matrix $C$ such that $M_k(h_k) =M_k(g_k) C$ and  $M_{k+1}(h_{k+1})= C^{-1}M_{k+1}(g_{k+1})$), we can use a complete ordering of the parameter space $\RR^{K\times S}$ and impose, using ${\mathcal M}$, the largest of all the equivalent versions of a parameters to be considered.

\subsection{Bibliography}
\subsubsection{Other matrix factorization and compressed sensing}
 The content of this paper is strongly related to and can be considered as an extension of the research field usually named {\em compressed sensing}. Because of the importance of this field of research and to simplify the reading for readers whose main interest is in deep learning, we have separated this part of the bibliography and placed it in \cref{bibli-ref}. Notice that the statement of \cref{unformal-stat} can be interpreted in the context of signal recovery. In particular, the results on deep structured linear networks can probably be specialized to be applicable to matrix factorization problems for which stability properties have not been established \cite{tsiligkaridis2013convergence,FTL_IJCV,FTO,magoarou2015flexible,LeMagoarouGraph,kondor2014multiresolution,Lyu:2013:ASM:2999611.2999679}. We have not investigated this potential. 

\subsubsection{Tensors and deep networks}
The analysis conducted in this paper is based on a connection, named {\em tensorial lifting}, between deep structured linear networks and a tensor problem (see \cref{lift-sec}). The tensorial lifting has already been described in \cite{MalgouyresLandsbergITW} but other connections between tensor and network problems have been described by other authors. In particular, in \cite{cohen2016expressive,cohen2016convolutional,khrulkov2017expressive}, the authors define a score function using a tensor. They highlight a network topology that computes the score function defined by a tensor decomposable using a CP-decomposition, a Hierarchical Tucker \cite{cohen2016expressive,cohen2016convolutional} or a tensor train decomposition\cite{khrulkov2017expressive}. They then deduce the expressive power of the network topology from the connections between the tensor decompositions. These results highlight and analyze why deep networks are more expressive than shallow ones. Tensors and tensor decomposition have also been used to represent the cross-moment and construct a solver \cite{janzamin2015beating}, encode the convolution layers with a tensor of order $4$ and manipulate this tensor to improve the network \cite{lebedev2014speeding,novikov2015tensorizing,yunpeng2017sharing}, to represent a tensor layer \cite{socher2013reasoning,yu2013deep}.

\subsubsection{Stability property for neural networks}

To the best of our knowledge, the articles establishing stability properties are 
\cite{arora2014provable,brutzkus2017globally,li2017convergence,sedghi2014provable,zhong2017recovery,MalgouyresLandsbergITW}.

Among them, \cite{brutzkus2017globally,li2017convergence,zhong2017recovery} consider a family of networks of depth $1$ or $2$ (depending on the article, the definition of the depth may vary). The article \cite{sedghi2014provable} contains a study on deep networks (the depth can be large), but the study only focuses on the recovery of one layer. The articles \cite{arora2014provable,MalgouyresLandsbergITW} consider networks without depth limitation.

In \cite{brutzkus2017globally}, the authors consider the minimization of the statistical risk (not the empirical risk). The input is assumed Gaussian and the output is generated by a network involving one linear layer followed by ReLU and a mean. The number of intermediate nodes is smaller than the input size. They provide conditions guaranteeing that, with high probability, a randomly initialized gradient descent algorithm converges to the true parameters. The authors of \cite{li2017convergence} consider a feed forward network made of one unknown linear layer, followed by ReLU and a sum. The size of the intermediate layer equals the size of the data, the size of the output is $1$. Again, they assume Gaussian input data and consider the minimization of the risk (not the empirical risk). They show that the stochastic gradient descent converges to the true solution. In \cite{zhong2017recovery}, the authors consider a non-linear layer followed by a linear layer. The size of the intermediate layer is smaller than the size of the input and the size of the output is $1$. They describe an initialization algorithm based on a tensor decomposition such that with high probability, the gradient algorithm minimizing the empirical risk converges to the true parameters that generated the data.

The authors of \cite{sedghi2014provable} consider a feed-forward neural network and show that, if the input is Gaussian or its distribution is known, a method based on moments and sparse dictionary learning can retrieve the parameers defining the first layer. Nothing is said about the stability or the estimation of the other layers.

The authors of \cite{arora2014provable} consider deep feed-forward networks which are very sparse and randomly generated. They show that they can be learned with high probability one layer after another. However, very sparse and randomly generated networks are not used in practice and one might want to study more versatile structures. The article \cite{MalgouyresLandsbergITW} studies deep structured linear networks (without the models $\Mod{}$) and uses the same tensorial lifting we use here. However, in \cite{MalgouyresLandsbergITW} the function $d$ measuring the error between parameters is only defined using the $\ell^\infty$ norm and ii not a metric. The transversality condition of \cite{MalgouyresLandsbergITW} is sufficient to guarantee the stability but is not necessary. All these weakness are corrected in this extended version. The general result is also specialized to deep convolutional linear networks.

\subsection{Organization of the paper}

Because it is strongly related, we give an extensive bibliography on compressive sensing and stable recovery properties for matrix factorization problems in \cref{bibli-ref}. We describe the framework of the paper and our notations in \cref{notation-sec}.

The main contributions of this paper are:
\begin{itemize}
\item In \cref{tensor-sec}, we investigate and recall several results on tensors, tensor rank and the Segre embedding. In particular, we investigate how the Segre embedding distorts distances.

\item In \cref{lift-sec}, we describe the {\em tensorial lifting}. It expresses any deep structured linear networks in a generic multilinear format. The latter composes a linear lifting operator and the Segre embedding.

\item When $\delta = \eta = 0$ (see \cref{ident-sec}):
\begin{itemize}
\item We establish a simple geometric condition on the intersection of two sets which is necessary and sufficient to guarantee the identifiability of the parameters up to scale ambiguity (\cref{ident-thm}).
\item  We provide simpler conditions which involve the rank of the lifting operator (defined in \cref{lift-sec}) such that:
\begin{itemize}
\item Under-specified case: If the lifting operator rank is large (e.g. larger than $2K(S-1)+2$, when $\Mod{} = \hS$) and the lifting operator is random, for almost every lifting operator, the solution of
\[M_1(h_1) \cdots M_K(h_K) = X \]
 is identifiable (\cref{suffic-ident-thm}).
\item Over-specified case: If the  lifting operator rank is small (e.g. smaller than $2S-1$, when $\Mod{} = \hS$), the solution of
\[M_1(h_1) \cdots M_K(h_K) = X \]
is not identifiable (\cref{necess-ident-thm});
\end{itemize}
\item We also provide a simple algorithm to compute the rank of the lifting operator (\cref{rkA-prop}).
\end{itemize}

\item Stability guarantee statements for deep structured linear networks are in \cref{stable-sec}:
\begin{itemize}
\item We define the \NSPlong~(Definition \cref{dnsp-def}): a generalization of the usual Null Space Property \cite{cohen2009compressed} that also applies to the deep problems.
\item We establish that the \NSPlong~ is a necessary and sufficient condition to guarantee stability (see the informal statement above or \cref{suf-stble-recovery} and \cref{nec-stble-recovery}).

\end{itemize}


\item We specialize the results to convolutional linear networks in \cref{conv-tree}  and establish a simple condition that can be computed (see Algorithm \cref{algo_check}). The is such that (see \cref{stabl-rec-network-thm})
\begin{itemize}

\item If the condition is satisfied the convolutional linear networks can be stably recovered;

\item If the condition is not satisfied, the convolutional linear network is not identifiable.
\end{itemize}
In simple words, the condition holds when the supports of the convolution kernels are sufficiently scattered. This is not satisfied by the convolutional kernel used in applications and explains their instability. 

\end{itemize}

\ifthenelse{\boolean{long}}
{}
{Because of space constraints, all the proofs are provided as a supplementary material or in the public archive \cite{MalgouyresLandsbergLong}.
}

\section{Bibliography on matrix factorization and compressed sensing}\label{bibli-ref}

Before describing the bibliography on compressed sensing, we interpret this stability statement of \cref{unformal-stat} in the context of signal processing. In signal processing, we  usually know that $\overline h$ exists and $\delta$ represents the sum of a modeling error and noise. The inequality \cref{erug} guarantees that, when the condition is satisfied, even an approximative minimizer of 
\begin{equation}\label{model}
 \argmin_{L\in\NN, (h_k)_{k=1..K} \in\Mod{L}} ~ \| M_1(h_1) \cdots M_K(h_K) - X \|^2.
\end{equation}
leads to a solution $h^*$ close $\overline h$. This property is often named: {\em stable recovery guarantee}.

When $\delta = 0$ (i.e., the data exactly fits the model and is not noisy) and $\eta=0$ (i.e., \cref{model} is perfectly solved) this is an {\em identifiability guarantee}. This is a necessary condition of stable recovery.

In this section, we distinguish the cases $K=1$, $K=2$ and $K\geq 3$. 
\subsection{ $K=1$: Linear inverse problems}

The simplest version consists of a model with
one layer (i.e., $K=1$) and $\Mod{}=\hS$. Recovering $h_1$ from $X$ is a linear inverse problem. The data $X$ can be vectorized to form a column vector and the operator $M_1$ simply multiplies the column vector $h_1$ by a fixed (rectangular) matrix. Typically, when the linear inverse problem is over-determined, the latter matrix has more rows than columns, the uniqueness of a solution to \cref{model} depends on the column rank of the matrix and the stable recovery constant depends on the smallest singular value of $M_1$.

When the matrix is not full colum rank, the identifiability and stable recovery for this problem has been intensively studied for many constraints $\Mod{}$. In particular, for sparsity constraints this is the compressed/compressive sensing problem (see the seminal articles \cite{candes2006robust,donoho2006compressed}). 
Some compressed sensing statements (especially the ones guaranteeing that any minimizer of the $\ell^0$ problem stably recovers the unknown problem) are special cases ($K=1$) of the statements provided in this paper. We will not give a complete review on compressed sensing but would like to highlight the Null Space Property described in \cite{cohen2009compressed}. The fundamental limits of compressed sensing (for a solution of the $\ell^0$ problem) have been analyzed in detail in \cite{bourrier2014fundamental}.

Although the main novelty of the paper is to investigate stable recovery properties for any $K\geq 1$, we specialize the statements made for $K\geq 1$ to the case $K=1$ in order to illustrate the new statements and to provide a way of comparison with well known results.

\subsection{$K=2$: Bilinear inverse problems and bilinear parameterizations}

When $k\geq 2$, the problem becomes non-linear because of the product in \cref{model}. This significantly complicates the analysis. What follows are the main instances studied in the literature when $K=2$.

\paragraph{Non-negative Matrix factorization (NMF) and low rank prior:}

In non-negative matrix factorization \cite{lee1999learning}, $M_1$ and $M_2$ map the entries in $h_1$ and $h_2$ at prescribed locations in the factors (say, one column after another). The constraints $\Mod{}{}$ imposes that all the entries in $h_1$ and $h_2$ are non-negative. The NMF has been widely used for many applications. 

Conditions guaranteeing that the factors provided by the NMF identify\footnote{Stable recovery is not established.} the correct factors (up to rescaling and permutation) were first established in the pioneering work \cite{donoho2003does}. To the best of our knowledge, this is the first paper addressing recovery guarantees for a problem of depth $K=2$. It emphasizes a separability condition that guarantees identifiability. The proof is purely geometric and relies on the analysis of inclusions of simplicial cones.
This result is significantly extended in \cite{laurberg2008theorems}. In this paper, the continuity of the NMF estimator is established. Concerning computational aspects, NMF is NP-complete \cite{vavasis2009complexity}. However, under the separability hypothesis of \cite{donoho2003does}, the solution of the NMF problem can be computed in polynomial time \cite{arora2012computing}.

We can slightly generalize\footnote{The interested readers can check that this generalization only leads to a small change of the Lifting operator introduced in \cref{lift-sec}. It is therefore done at no cost.} the problem and introduce a linear degradation operator
\[H: \RR^{m\times n} \longrightarrow \RR^{m \times n}.
\]
Use the same mapping $M_1$ and $M_2$ as for the NMF, with $\Mod{}=\hS$, but with a small number of lines (resp columns) in $M_2(h_2)$ (resp. $M_1(h_1)$). Any solution of the problem
\[(h_1^*,h_2^*)\in \argmin_{(h_1,h_2)\in\hS} ~ \| H(M_1(h_1)M_2(h_2)) - X \|^2
\]
leads to a low rank approximation $M_1(h_1^*)M_2(h_2^*)$ of an inverse of $H$, at $X$. Again, a large corpus of literature exists on the low rank prior \cite{recht2010guaranteed,candes2009exact,fazel2008compressed,candes2010power}. 


\paragraph{Phase Retrieval:}
  Phase retrieval fits the framework described in the present paper when we take
\[M_1(h_1) = \DIAG{\calF h_1}\qquad M_2(h_2)= (\calF h_2)^*
\]
and
\[\Mod{}{} = \{(h,h)\in\hS\mid h\in\RR^S \}
\]
where $S$ is the size of the signal, $\calF$ computes $N$ linear measurements of any element in $\RR^S$ (typically Fourier measurements), $\DIAG{.}$ creates an $N\times N$ diagonal matrix whose diagonal contains the input and $^*$ is the (entry-wise) complex conjugate.

The tensorial lifting at the core of the present paper generalizes the lifting used in the inspiring work on PhaseLift \cite{li2013sparse,candes2013phaselift,candes2015phase}. As is often the case when $K=2$, PhaseLift is a semidefinite program that can be efficiently solved when the unknown is of moderate size. These papers also provide conditions on the measurements guaranteeing that the phases are stably recovered by PhaseLift.

The benefit of the generalization introduced with the tensorial lifting is that it applies to any multilinear inverse problem.

\paragraph{Self-calibration and de-mixing}

Measuring operators often depend linearly on parameters that are not perfectly known. The estimation of these parameters is crucial to restore the data measured by the device. This is the self-calibration problem. This naturally fits the setting of this article: - we let $h_1$ be the parameters defining the sensing matrix and $M_1(h_1)$ be the sensing matrix. Then $h_2$ defines the signal (or signals) contained in the column(s) of $M_2(h_2)$. 

Many instances of this problem have been studied and much progress has been made to obtain algorithms that can be applied to problems of larger and larger size. This leads to a very interesting line of research. 

To the best of our knowledge, the first stable recovery statements concern the blind-deconvolution problem. In \cite{ahmed2014blind}, the authors use a lifting to transform the blind-deconvolution problem into a semidefinite program with an unknown whose size is the product of the sizes\footnote{With our notations this is simply $S\times S$ but this can be much more favorable.} of $h_1$ and $h_2$. Such problems can be solved for unknowns of moderate size. The authors of \cite{ahmed2014blind} provide explicit conditions guaranteeing the stable recovery with high probability. This idea has been generalized and applied to other similar problems in \cite{choudhary2014identifiability,bahmani2015lifting}. The authors of \cite{ling2015self} consider a significantly more general calibration model. In this model, $M_1(h_1)$ is diagonal and its diagonal contains the entries of $h_1$. $M_2(h_2)$ simply multiplies $h_2$ by a fixed known matrix (the theorems consider a random matrix). The constraint imposes $h_2$ to be sparse. For this problem, they prove that with high probability  the numerical method called SparseLift is stable with a controlled accuracy. SparseLift returns the left and right singular vectors of the solutions of an $\ell^1$ optimization problem whose unknown is the same as in \cite{ahmed2014blind}. However, solving an $\ell^1$ minimization problem is much simpler than a semi-definite problem. This is a very significant practical improvement.

As emphasized in \cite{DBLP:journals/corr/LiLSW16}, in order to motivate its non-convex approach, the only drawback of the numerical methods described in \cite{ahmed2014blind,ling2015self} is their complexity. The extra complexity is due to the fact that they optimize a variable in the product space $\RR^{S\times S}$ and then deduce an approximate solution of the un-lifted problem. This is what motivates the authors of \cite{DBLP:journals/corr/LiLSW16} to propose a non-convex approach. The constructed algorithm provably stably recovers the sensing parameters and the signals  with a geometric onvergence rate.

\paragraph{Sparse coding and dictionary learning:}

Sparse coding and dictionary learning is another kind of bilinear problem (see \cite{RubinsteinBrucksteinElad} for an overview). In that framework, the columns of $X$ contain the data. Most often, people consider two layers: $K=2$. The layer $M_1(h_1)$ is an optimized dictionary of atoms defined by the parameters $h_1$ and each column of $M_2(h_2)$ contains the code (or coordinates) of the corresponding column in $X$. Most often, $h_2$ is assumed sparse. 

The identifiability and stable recovery of the factors has been studied in many dictionary learning contexts and provides guarantees on the approximate recovery of both an incoherent dictionary and sparse coefficients  when the number of samples is sufficiently large (i.e., in our setting when $n$ is large).  In \cite{GribonvalSchnass}, the authors developed local optimality conditions in the noiseless case, as well as sample complexity bounds for local recovery when $M_1(h_1)$ is square and $M_2(h_2)$ are iid Bernoulli-Gaussian. This was extended to overcomplete dictionaries in \cite{Geng} (see also \cite{Schnass} for tight frames) and to the noisy case in \cite{Jenatton}.
 The authors of \cite{Spielman} provide exact recovery results for dictionary learning, when the coefficient matrix has Bernoulli-Gaussian entries and the dictionary matrix has full column rank. This was extended to overcomplete dictionaries in \cite{Agarwal} and in \cite{Arora} but only for approximate recovery. Finally, \cite{7088631} provides such guarantees under general conditions which cover many practical settings.

\paragraph{Contributions in these frameworks}
The present article considers the identifiability and stability of the recovery for any $K\geq 1$ in a general and unifying framework. As was already mentioned, we do not investigate computational issues. As will appear, later  the paper, the analogue of the lifting at the core of the algorithms described in the above papers (in particular the papers on phase retrieval and self-calibration) is a {\em tensorial lifting} (see \cref{lift-sec}) and involves tensors that cannot be manipulated in practice. Even when we are able to manipulate the tensors, the computation of the best rank $1$ approximation of such tensors is an open non-convex problem. Therefore, there is no numerically efficient and reliable way to extract the un-lifted parameters from an optimized tensor. Because of that, we have not yet pursued the construction of a numerical scheme based on the tensorial lifting when $K\geq 3$. As was already mentioned, at this writing, the success of algorithms for $K\geq 3$ is mostly supported by empirical evidence. Proving their efficiency is a wide open problem (see \cite{livni2014computational,haeffele2015global,choromanska2015loss,choromanska2015open,safran2016quality,baldi1989neural,baldi1995learning,kawaguchi2016deep,2018arXiv180206384V}). The purpose of the paper is to provide guarantees on the stability of the solution when such an empirical success occurs.

The specialization of the presented results to problems with $K=2$ leads to necessary and sufficient conditions for the stable recovery. This is slightly different from the usual approach. Usually, authors provide sufficient conditions and argue their sharpness by comparing the number of samples required by their method and the information theoretic limit (typically, the number of independent variables of the problem).

It would of course be interesting to see how far it is possible to unify the different problems with $K=2$ using the framework of this paper. We have however not pursued this route and instead focus on the situation $K\geq 3$. 

\subsection{$K\geq 3$}
The difficulties,  when $K\geq 3$, come from the fact that tools used for problems with $K=2$ are not applicable. In particular, we cannot use the usual lifting, the singular value decomposition or the sin-$\theta$ theorem in \cite{davis1970rotation}. Often, these tools are replaced by analogous objects involving tensors. This complicates the analysis and prohibits the use of numerical schemes that manipulate lifted variables.
 
To the best of our knowledge, little is known concerning the identifiability and the stability of  matrix factorization when $K\geq 3$. The uniqueness of the factorization corresponding to the Fast Fourier Transform was proved in \cite{LeMagoarouThesis}. Other results  consider the identifiability of the factors which are sparse and random \cite{neyshabur2013sparse}. The authors of the present paper have announced preliminary versions of the results described here in \cite{MalgouyresLandsbergITW}. They are significantly extended here.


\section{Notation and summary of the hypotheses}\label{notation-sec}

We continue to use   the notation introduced in the introduction.
For an integer $k\in\NN$, set $\NNN{k} = \{1,\cdots, k\}$.

We consider $K\geq 1$ and $S\geq 2$ and real-valued tensors of order $K$ whose axes are of size $S$, denoted $T\in\RR^{S\times \cdots \times S}$. The space of tensors is
 abbreviated  $\TS$. The entries of $T$ are denoted $T_{i_1,\cdots,i_K}$, where $(i_1,\cdots, i_K)\in\iS$.
 For $\ibf\in \iS$, the entries of $\ibf$   are $\ibf= (i_1,\cdots, i_K)$ (for $\jbf\in \iS$ we let  $\jbf= (j_1,\cdots, j_K)$, {\it etc}). We
 either write  $T_\ibf$ or $T_{i_1,\cdots,i_K}$.

To simplify notations, from now on, the parameters defining the factors are gathered in a single matrix and are denoted with bold fonts $\hbf\in\hS$. The $k^{\mbox{th}}$ vector containing the parameters for the layer $k$ is denoted   $\hbf_k\in\RR^S$. The $i^{\mbox{th}}$ entry of the $k^{\mbox{th}}$ vector is denoted   $\hbf_{k,i}\in\RR$. A vector not related to an element in $\hS$ is denoted $h\in\RR^S$ (i.e., using a light font). Throughout the paper we assume
\[\Mod{} = (\Mod{L})_{L\in\NN}\mbox{, with }\Mod{L}\subset \hS.
\]
We also assume that, for all $L\in\NN$, $\Mod{L}\neq\emptyset$. They can however be equal or constant after a given $L'$.

All the vector spaces $\TS$, $\hS$, $\RR^S$ {\it etc.} are equipped with the usual Euclidean norm. This norm is denoted $\|.\|$ and the scalar product $\PS{.}{.}$. In the particular case of matrices, $\|.\|$   corresponds to the Frobenius norm.  We also use the usual $p$ norm, for $p\in[1,\infty]$, and denote it by $\|.\|_p$. In particular, for $\hbf\in\hS$ and $T\in \TS$, we have for $p<+\infty$
\[\|\hbf\|_p =  \left( \sum_{k=1}^K\sum_{i=1}^S |\hbf_{k,i}|^p \right)^{1/p} \qquad\mbox{,}\qquad \|T\|_p =  \left( \sum_{\ibf\in\iS} |T_\ibf|^p \right)^{1/p} \]
and
\[\|\hbf\|_{\infty} =  \max_{ \substack{k\in\NNN{K} \\ i\in\NNN{S}} } |\hbf_{k,i}| \qquad\mbox{,}\qquad \|T\|_{\infty} =   \max_{\ibf\in\iS} |T_\ibf|. \]

Set
\begin{equation}\label{eruonoqen}
\hS_* = \{\hbf \in \hS\mid \forall k\in\NNN{K}, \|\hbf_k\| \neq 0  \}.
\end{equation}
Define an  equivalence relation on $\hS_*$: for any  $\hbf$,  ${\mathbf g} \in \hS$, $\hbf \sim {\mathbf g}$  if and only if there exist $(\lambda_k)_{k\in\NNN{K}} \in\RR^K$ such that
\begin{equation}\label{equiv-rel}
\prod_{k=1}^K \lambda_k = 1 \qquad \mbox{ and }  \qquad \hbf_k = \lambda_k {\mathbf g}_k, \forall k\in\NNN{K}.
\end{equation}
Denote the equivalence class of  $\hbf\in\hS_*$ by $\CL{\hbf}$.

The zero tensor is of rank $0$. A non-zero tensor $T\in\TS$ is  of {\it rank $1$} (or decomposable) if and only if there exists  $\hbf\in\hS_*$ such that $T$ is the outer product of the vectors $\hbf_k$, for $k\in\NNN{K}$.
That is,  for any $\ibf\in \iS$,
\[T_\ibf = \hbf_{1,i_1} \cdots \hbf_{K,i_K}.
\]
Let $\Sigma_1\subset \TS$ denote the set of tensors of rank $0$ or  $1$.

The {\it rank}  of a tensor $T\in\TS$ is
\[\RK{T}=\min\{ r \in\NN \mid \mbox{ there exists } T_1,\cdots, T_r\in\Sigma_1
\mbox{ such that }T=T_1+\cdots+ T_r  \}.
\]

For   $r\in\NN$, let
\[\Sigma_r = \{T\in\TS \mid  \RK{T} \leq r\}.
\]



The $*$ \underline{super}script refers to optimal solutions. A set with a $*$ \underline{sub}script means that $0$ is ruled out of the set. In particular, $\Sigma_{1,*}$ denotes the non-zero tensors of rank $1$. Attention should   be paid to $\hS_*$  (see \cref{eruonoqen}).


\section{Facts on the Segre embedding and tensors of rank $1$ and $2$} \label{tensor-sec}

Parametrize  $\Sigma_1\subset \TS$ by   the map
\begin{equation}\label{defP}
\begin{array}{rcl}
P: \hS & \longrightarrow & \Sigma_1 \subset \TS \\
\hbf & \longmapsto & (\hbf_{1,i_1}\hbf_{2,i_2} \cdots \hbf_{K,i_K})_{\ibf\in\iS}.
\end{array}
\end{equation}

The map $P$ is called the Segre embedding and is often denoted $\widehat{Seg}$ in the algebraic geometry literature.

{\bf Standard Facts:}
\begin{enumerate}
\item \label{rk1-Ident-cor}{\bf Identifiability of $\CL{\hbf}$ from $P(\hbf)$:}
For   $\hbf$ and $\gbf\in\hS_*$,  $P(\hbf)=P(\gbf)$ if and only if $\CL{\hbf} = \CL{\gbf}$.
\item  \label{rk1-geo-thm}{\bf Geometrical description of $\Sigma_{1,*}$: }
$\Sigma_{1,*}$ is a smooth (i.e., $C^\infty$) manifold of dimension $K(S-1)+1$ (see, e.g.,  \cite{landsberg2012tensors}, chapter 4, pp. 103).

\item \label{rk2-geo-thm}{\bf Geometrical description of $\Sigma_{2}$:}
We recall that the singular locus $(\overline{\Sigma}_2)_{sing}$ of the closure $\overline{\Sigma}_2$ of $\Sigma_2$ has dimension strictly less than that of $\overline{\Sigma}_2$ and that $\overline{\Sigma}_2\backslash (\overline{\Sigma}_2)_{sing}$ is a smooth manifold. The dimension of $\overline{\Sigma}_2\backslash (\overline{\Sigma}_2)_{sing}$ is $2K(S-1)+2$  when $K>2$,  and   is $4(S-1)$ when $K=2$
(see, e.g.,  \cite{landsberg2012tensors}, chapter 5).
\end{enumerate}


We can improve Standard Fact \cref{rk1-Ident-cor} and obtain a stability result guaranteeing that if we know a rank $1$ tensor sufficiently close to $P(\hbf)$, we  approximately know $\CL{\hbf}$. In order to state this, we need to define a metric on $\hS_* /\sim$ (where $\sim$ is defined by \cref{equiv-rel}). This has to be considered with care since, whatever $\hbf \in \hS_*$, the subset
$\{ h\mid h\in \CL{\hbf}\} $ is not compact. In particular, considering
$$\hbf'_k = \left\{\begin{array}{ll}
\lambda~ \hbf_k & \mbox{ if } k=1 \\
\lambda^{-\frac{1}{K-1}}~ \hbf_k& \mbox{ otherwise}
\end{array}\right.
$$
when $\lambda$ goes to infinity, we easily construct examples that make
 the standard metric on  equivalence classes useless\footnote{For instance, if $\hbf$ and $\gbf\in\hS_*$ are such that  $\hbf_1=\gbf_1$, we have $$\inf_{\hbf'\in\CL{\hbf},\gbf'\in\CL{\gbf}  } \|\hbf'-\gbf'\|_p=0$$
even though we might have $\hbf_2\neq\gbf_2$ (and therefore $\CL{\hbf} \neq \CL{\gbf}$). This does not define a metric.

Also, when $\hbf$ and $\gbf$ are such that $\hbf_k\neq\gbf_k$, whatever $k\in\NNN{K}$, we have
$$\sup_{\hbf'\in\CL{\hbf}} \inf_{\gbf'\in\CL{\gbf}  } \|\hbf'-\gbf'\|_p=+\infty.$$
Therefore, the Hausdorff distance between $\CL{\hbf}$ and $\CL{\gbf}$ is infinite for almost every pair $(\hbf,\gbf)$. This metric is therefore not very useful in the present context.}.

This leads us to consider
\[\hSdiag = \{\hbf \in \hS_* \mid \forall k\in\NNN{K}, \|\hbf_k\|_\infty = \|\hbf_1\|_\infty  \}.\]
The interest in  this set comes from the fact that, whatever $\hbf \in \hS_*$, the set $\CL{\hbf}\cap\hSdiag$ is finite. Indeed, if $\gbf\in \CL{\hbf}\cap\hSdiag$ the $(\lambda_k)_{k\in\NNN{K}} \in\RR^K$ such that, for all $k\in\NNN{K}$, $\hbf_k = \lambda_k \gbf_k$ must all satisfy $|\lambda_k|=1$, i.e., $\lambda_k=\pm1$.

\begin{defi}\label{def_d}
For any $p\in[1,\infty]$, we define the mapping   $d_p:(\hS_* /\sim \times \hS_* /\sim) \rightarrow \RR$ by
\[d_p(\CL{\hbf}, \CL{\gbf}) =  \inf_{\substack{\hbf'\in\CL{\hbf}\cap\hSdiag \\
\gbf'\in\CL{\gbf} \cap\hSdiag  }} \|\hbf'-\gbf'\|_p\qquad\forall \hbf,\, \gbf \in \hS_*.
\]
\end{defi}

\begin{prop}\label{metric-prop}
For any $p\in[1,\infty]$, $d_p$ is a metric on $\hS_* /\sim$.
\end{prop}
\ifthenelse{\boolean{long}}
{The proof is in \cref{metric-prop-proof}.}
{The proof is in the supplementary material \cref{metric-prop-proof} and the public archive \cite{MalgouyresLandsbergLong}.
}

Notice that the equivalence relationship and metric defined above are not adapted to operators $M_k$ allowing invariance such as permutations. More precisely, for some operators $M_k$, there exists $\hbf$, $\gbf$ and a permutation matrix $C$ such that $(\hbf_k,\hbf_{k+1})\neq (\gbf_k,\gbf_{k+1})$ and  $M_k(\hbf_k) =M_k(\gbf_k) C$ and  $M_{k+1}(\hbf_{k+1})= C^{-1}M_{k+1}(\gbf_{k+1})$. In such a case, we have $d_p(\CL{\hbf}, \CL{\gbf}) \neq 0$. However, the features defined at the layer $k$ are just permuted and can still be interpreted. As already said, in such a case, it is possible to use the models $\Mod{}$ to select one the equally interpretable $\hbf$.

Using the above metric, we can state that not only $\CL{\hbf}$ is uniquely determined by $P(\hbf)$, but this operation is stable.

\begin{thm}\label{rk1-Ident-thm}{\bf Stability of $\CL{\hbf}$ from $P(\hbf)$}

Let $\hbf$ and $\gbf\in\hS_*$ be such that
$\|P(\gbf)-P(\hbf)\|_\infty\leq \frac{1}{2}  ~ \max\left(\|P(\hbf)\|_\infty,\|P(\gbf)\|_\infty\right) $. For all $p, q\in[1,\infty]$,
\begin{equation}\label{Pinv-lipschitz}
d_p(\CL{\hbf}, \CL{\gbf}) \leq 7 (KS)^{\frac{1}{p}} \min\left(\|P(\hbf)\|_{\infty}^{\frac{1}{K}-1},\|P(\gbf)\|_{\infty}^{\frac{1}{K}-1} \right) \|P(\hbf)-P(\gbf)\|_{q}.
\end{equation}
\end{thm}
\ifthenelse{\boolean{long}}
{The proof of the theorem is in \cref{rk1-Ident-thm-proof}.}
{The proof of the theorem is in the supplementary material \cref{rk1-Ident-thm-proof} and in \cite{MalgouyresLandsbergLong}.
}

In the final result, the bound established in \cref{rk1-Ident-thm} plays a role similar to the $sin-\theta$ Theorem of \cite{davis1970rotation} in \cite{ling2015self,candes2013phaselift,ahmed2014blind}.


The following proposition shows that the upper bound in \cref{Pinv-lipschitz} cannot be improved by a  significant factor, in particular when $q$ is large.

\begin{prop}\label{rk1-Ident-sharp}
There exist  $\hbf$ and $\gbf\in\hS_*$ such that $\|P(\gbf)\|_\infty \leq \|P(\hbf)\|_\infty$, $\|P(\gbf)-P(\hbf)\|_\infty\leq \frac{1}{2}  ~ \|P(\hbf)\|_\infty$ and
\[7 (KS)^{\frac{1}{p}} \|P(\hbf)\|_{\infty}^{\frac{1}{K}-1} \|P(\hbf)-P(\gbf)\|_{q}\leq C_q ~ d_p(\CL{\hbf}, \CL{\gbf}),
\]
where
\[C_q = \left\{\begin{array}{ll}
28 (KS)^{\frac{1}{q}} & \mbox{ if } q <+\infty ,\\
28& \mbox{ if } q=+\infty.
\end{array}\right.
\]
\end{prop}
\ifthenelse{\boolean{long}}
{The proof of the theorem is in \cref{rk1-Ident-sharp-proof}.}
{The proof of the proposition is in the supplementary material cref{rk1-Ident-sharp-proof} and in \cite{MalgouyresLandsbergLong}.
}

As stated in the following theorem, we have a  more valuable  upper bound in the general case.

\begin{thm}\label{PLip-thm}{\bf \lq\lq Lipschitz continuity\rq\rq of $P$}

For any $q\in[1,\infty]$ and any $\hbf$ and $\gbf\in\hS_*$,
\begin{equation}\label{P-lipschitz}
 \|P(\hbf)-P(\gbf)\|_{q} \leq S^{\frac{K-1}{q}} K^{1-\frac{1}{q}}  \max\left(\|P(\hbf)\|_{\infty}^{1-\frac{1}{K}} , \|P(\gbf)\|_{\infty}^{1-\frac{1}{K}} \right) d_q(\CL{\hbf}, \CL{\gbf}).
\end{equation}
\end{thm}
\ifthenelse{\boolean{long}}
{The theorem is proved in \cref{PLip-thm-proof}.}
{The theorem is proved in the supplementary material \cref{PLip-thm-proof} and in \cite{MalgouyresLandsbergLong}.
}

Notice that, considering $\hbf$ and $\gbf\in\hS$ such that $\hbf_{k,i} = 1$ and $\gbf _{k,i} = \varepsilon$, for all $k\in\NNN{K}$ and $i\in\NNN{S}$ and for a $0<\varepsilon\ll 1$, we easily calculate
\[S^{\frac{K-1}{q}} K^{1-\frac{1}{q}}  \max\left(\|P(\hbf)\|_{\infty}^{1-\frac{1}{K}} , \|P(\gbf)\|_{\infty}^{1-\frac{1}{K}} \right) d_q(\CL{\hbf}, \CL{\gbf}) \leq K  \|P(\hbf)-P(\gbf)\|_{q}.
\]
As a consequence, the upper bound in \cref{PLip-thm} is tight up to at most a factor $K$.


\section{The tensorial lifting }\label{lift-sec}


The following proposition is clear (it can be shown by induction on $K$):

\begin{prop}
Let $M_k$, $k\in\NNN{K}$ be as in \cref{mkdef}. The entries  of the matrix $$M_1(\hbf_1)M_2(\hbf_2) \cdots M_K(\hbf_K)$$ are multivariate polynomials whose variables are the entries of $\hbf \in\hS$. Moreover, every
entry is the sum of monomials of degree $K$. Each monomial is a constant times  $\hbf_{1,i_1} \cdots \hbf_{K,i_K}$, for some $\ibf\in \iS$.
\end{prop}

Notice that any monomial $\hbf_{1,i_1} \cdots \hbf_{K,i_K}$ is the entry $P(\hbf)_\ibf$ in the tensor $P(\hbf)$. Therefore every polynomial in the previous proposition takes the form $\sum_{\ibf\in\iS} c_\ibf P(\hbf)_\ibf$ for some constants $(c_\ibf)_{\ibf\in\iS}$ independent of $\hbf$. In words, every entry  of the matrix $M_1(\hbf_1)M_2(\hbf_2) \cdots M_K(\hbf_K)$ is obtained by applying a linear form to $P(\hbf)$. Moreover, the polynomial coefficients defining the linear form are uniquely determined by the linear maps $M_1$, $\cdots$, $M_K$. This leads to the following statement.
\begin{cor}{\bf Tensorial Lifting} \label{AP=M-cor}

Let $M_k$, $k\in\NNN{K}$ be as in \cref{mkdef}. The map
$$(\hbf_1,\hdots, \hbf_K)\longmapsto  M_1(\hbf_1)M_2(\hbf_2) \cdots M_K(\hbf_K),
$$
uniquely determines a linear map
\[\calA:\TS \longrightarrow \RR^{m\times n},
\]
such that for all $\hbf\in\hS$
\begin{equation}\label{AP=M}
M_1(\hbf_1)M_2(\hbf_2) \cdots M_K(\hbf_K) = \calA P(\hbf).
\end{equation}
\end{cor}

We call \cref{AP=M} and its use the {\em tensorial lifting}. When $K=1$, we simply have $\calA = M_1$. When $K=2$ it corresponds to the usual lifting already exploited to establish stability results for phase recovery, blind-deconvolution, self-calibration, sparse coding, {\it etc}. Notice that, when $K\geq 2$, it may be difficult to provide a closed form expression for the operator $\calA$. We can however
determine simple properties of  $\calA$. In most reasonable cases, $\calA$ is sparse. If the operators $M_k$ simply embed the values of $h$ in a matrix, the matrix representing $\calA$ only contains zeros and ones. Since the operators $M_k$ are known, we can compute $\calA P(\hbf)$, for any $\hbf\in\hS$, using \cref{AP=M}. Said differently, we can compute $\calA$ for any rank $1$ tensor.
Therefore, since $\calA$ is linear, we can compute $\calA T$ for any low rank tensor $T$. If the dimensions of the problem permit, one can  manipulate $\calA$ in a basis of $\TS$.

Since $\RK{\calA}$ is an important quantity, we emphasize that  $\RK{\calA}\leq mn$. It is also possible to compute $\RK{\calA}$, when $m  n$  is not too large, using the following proposition.
\begin{prop}\label{rkA-prop}
For  $R$ independent random $\hbf^r$, with $r=1..R$, according to the normal distribution in $\hS$, we have with probability $1$
\begin{equation}\label{ejvn}
\dim(\SPAN{ (\calA P(\hbf^r))_{r = 1..R}}) =\left\{\begin{array}{ll}
R & \mbox{ if } R\leq \RK{\calA} \\
\RK{\calA} & \mbox{ otherwise.}
\end{array}\right.
\end{equation}
\end{prop}

\ifthenelse{\boolean{long}}
{The proof is in \cref{rkA-prop-proof}.}
{The proof is in the supplementary material \cref{rkA-prop-proof} and the public archive \cite{MalgouyresLandsbergLong}.
}

Using Corollary \cref{AP=M-cor}, when \cref{model} has a minimizer, we rewrite in the form
\begin{equation}\label{model1}
\hbf^* \in \argmin_{L\in\NN, \hbf\in\Mod{L}} ~ \|\calA P(\hbf) - X \|^2.
\end{equation}

We now decompose this problem into two sub-problems: A least-squares problem
\begin{equation}\label{Pb1}
T^* \in \argmin_{T\in\TS}  \|\calA T - X \|^2
\end{equation}
and a non-convex problem
\begin{equation}\label{Pb2}
\hbf'^* \in \argmin_{L\in\NN, \hbf\in\Mod{L}}  \|\calA (P(\hbf) -T^*) \|^2.
\end{equation}

\begin{prop}\label{lifting-prop}
Let $X$ and $\calA$ be such that \cref{model} has a minimizer:
\begin{enumerate}
\item Let $\hbf^*$ be  a  solution of \cref{model1}. Then, for any solution $T^*$ of \cref{Pb1}, $\hbf^*$ also minimizes \cref{Pb2}.
\item Let $T^*$ be  a solution of \cref{Pb1} and $\hbf'^*$ a  solution of \cref{Pb2}. Then, $\hbf'^*$ also minimizes \cref{model1}.
\end{enumerate}
\end{prop}

\ifthenelse{\boolean{long}}
{The proof is in \cref{lifting-Prop-proof}.}
{The proof is in the supplementary material \cref{lifting-Prop-proof} and the public archive \cite{MalgouyresLandsbergLong}.
}

From now on, because of the equivalence between solutions of \cref{Pb2} and \cref{model1}, we stop using the notation $\hbf'^*$ and write $\hbf^*\in\argmin_{L\in\NN, \hbf\in\Mod{L}}  \|\calA (P(\hbf) -T^*) \|^2 $.

\section{Identifiability (error free case)}\label{ident-sec}

Throughout this section, we assume that $X$ is such that there exists $\overline{L}$ and $\overline{\hbf}\in\Mod{\overline{L}}$ such that
\begin{equation}\label{noiseless}
X=  M_1(\overline{\hbf}_1) \cdots M_K(\overline{\hbf}_K).
\end{equation}
Under this assumption, $X=\calA P(\overline{\hbf})$, so
\[P(\overline{\hbf})\in \argmin_{T\in\TS}  \|\calA T - X \|^2.
\]
 Moreover, we trivially have $P(\overline{\hbf})\in\Sigma_1$ and therefore $\overline{\hbf}$ minimizes \cref{Pb2}, \cref{model} and \cref{model1}. As a consequence, \cref{model} has a minimizer.

We ask   whether there exist  guarantees that the resolution of \cref{model}
allows one  to recover $\overline{\hbf}$ up to the usual uncertainties.

In this regard, for any $\hbf\in\CL{\overline{\hbf}}$, we have $P(\hbf) = P(\overline{\hbf})$ and therefore $\calA P(\hbf) = \calA P(\overline{\hbf})=X$.
Thus unless we make further assumptions on $\overline{\hbf}$, we cannot expect to distinguish any particular element of $\CL{\overline{\hbf}}$ using only $X$. In other words, recovering   $\CL{\overline{\hbf}}$ is the best we can hope for.

\begin{defi}{\bf Identifiability}

We say that {\em $\CL{\overline{\hbf}}$ is identifiable} if  the elements of $\CL{\overline{\hbf}}$ are the only    solutions of \cref{model}.

We say that {\em $\Mod{}$ is identifiable} if for every $L\in\NN$ and every $\overline{\hbf}\in\Mod{L}$, $\CL{\overline{\hbf}}$ is identifiable.
\end{defi}

\begin{prop}{\bf Characterization of the global minimizers}\label{carac-thm}

For any $L^*\in\NN$ and any $\hbf^*\in\Mod{L^*}$, $(L^*,\hbf^*)\in\argmin_{L\in\NN, \hbf\in\Mod{}}  \|\calA P(\hbf) -X \|^2 $ if and only if
\[P({\hbf}^*) \in P(\overline{\hbf}) +\KER{\calA}.
\]
\end{prop}
\ifthenelse{\boolean{long}}
{The Proposition is proved in \cref{carac-thm-proof}.}
{The Proposition is proved in the supplementary material \cref{carac-thm-proof} and in \cite{MalgouyresLandsbergLong}.
}

In order to state the following proposition, we define for any $L$ and $L'\in\NN$
\[P(\Mod{L}) - P(\Mod{L'}):=  \left\{ P(\hbf) - P(\gbf) \mid \mbox{  \ } \hbf\in\Mod{L} \mbox{ and } \gbf \in\Mod{L'} \right\} \subset \TS.
\]

\begin{prop}{\bf Necessary and sufficient conditions of identifiability}\label{ident-thm}

\begin{enumerate}
\item For any $\overline{L}$ and $\overline{\hbf}\in\Mod{\overline{L}}$: $\CL{\overline{\hbf}}$ is identifiable if and only if for any $L\in\NN$
\[\bigl(P(\overline{\hbf})+\KER{\calA}\bigr) \cap P(\Mod{L}) ~\subset ~\{ P(\overline{\hbf})\}.
\]
\item $\Mod{}$ is identifiable if and only if for any $L$ and $L'\in\NN$
\begin{equation}\label{ident-hS}
\KER{\calA} \cap \bigl( P(\Mod{L}) - P(\Mod{L'}) \bigr)  \subset \{0\}.
\end{equation}
\end{enumerate}

\end{prop}
\ifthenelse{\boolean{long}}
{The proposition is proved in \cref{ident-thm-proof}.}
{The proposition is proved in the supplementary material \cref{ident-thm-proof} and in \cite{MalgouyresLandsbergLong}.
}

In the context of the usual compressed sensing (i.e., when $K=1$, $\Mod{}$ contains $L$-sparse signals, $\calA$ is a rectangular matrix with full row rank and $X$ is a vector), the proposition is already stated in Lemma 3.1 of \cite{cohen2009compressed}.

In reasonably small cases and when $P(\Mod{})$ is algebraic, one can  use tools from numerical algebraic geometry such as those described in \cite{Hauenstein20103349,Hauenstein20136809} to check whether the condition \cref{ident-hS} holds or not. The drawback of \cref{ident-thm} is that, given a deep structured linear network as described by $\calA$, the condition \cref{ident-hS} might be difficult  to verify.

We therefore establish simpler conditions related to the identifiability of $\Mod{}$. First we establish a condition such that for almost every $\calA$ satisfying it, $\Mod{}$ is identifiable. The main benefit of this condition is that its constituents  can be computed in many practical situations.

Before that, we recall a few facts of algebraic geometry, for $X,Y\subset \RR^N$, the {\it join} of $X$ and $Y$ (see, e.g., \cite[Ex. 8.1]{Harris})
is
$$J(X,Y):=\overline{\{ sx+ty \mid x\in X,\ y\in Y, \ s,t\in \RR\}}.
$$
If for all $L\in\NN$, $\Mod{L}$   is Zariski closed and invariant under rescaling
(e.g., if they are all linear spaces),
then $P(\Mod{L}) - P(\Mod{L'})$ is
a Zariski open subset of $J(P(\Mod{L}),  P(\Mod{L'}))$. In general, it is contained
in this join.

Recall the following fact (*):   for complex algebraic varieties  $X,Y\subset \CC^N$,
   any component $Z$  of $X\cap Y$ has $\DIM Z\geq \DIM X+\DIM Y-N$,
 and   equality holds generically (we make   \lq\lq generically\rq\rq\ precise
 in our context below).
 Moreover, if $X,Y$ are invariant under rescaling, since $0\in X\cap Y$,
 we have $X\cap Y\neq \emptyset$. (See, e.g., \cite[\S I.6.2]{MR3100243}.)

This intersection result indicates that if there exists $L,L'$ such that 
 $$\rk(\calA) < \DIM {P(\Mod{L}) - P(\Mod{L'})}$$ we expect to have non-identifiability; and if the rank is larger, for all pair $L,L'$, we expect identifiability. More precisely:

\begin{thm}{\bf Almost surely sufficient condition for Identifiability}\label{suffic-ident-thm}

For almost every $\calA$ such that  
$$\rk(\calA) \geq \DIM {J(P(\Mod{L}) , P(\Mod{L'}))},\qquad \mbox{ for all }L, L',
$$  
$\Mod{}$ is identifiable. 
\end{thm}
\ifthenelse{\boolean{long}}
{The theorem is proved in \cref{suffic-ident-thm-proof}.}
{The theorem is proved in the supplementary material \cref{suffic-ident-thm-proof} and in \cite{MalgouyresLandsbergLong}.
}

Since  $\DIM {J(P(\Mod{L}), P(\Mod{L'}))}\leq \DIM {P(\Mod{L})} +\DIM {P(\Mod{L'})}$,
if $D_{max}$ is the maximum dimension of $P(\Mod{L})$ over all $L$,
one has the same conclusion if $\rk(\calA) \geq 2D_{max}$.

When $K=1$, we  illustrate this result by interpreting it in
 the context of   compressive sensing,  where   $\hbf$ is a vector, $X$ is a vector, $\calA$ is a rectangular sampling matrix of full row rank and $\KER{\calA}$ is large. The statement  analogous  to   \cref{suffic-ident-thm}  in the compressive sensing framework takes the form: \lq\lq For almost every sampling matrix, any $L$ sparse signal $\hbf$ can be recovered from $\calA\hbf$ as soon as $2L\leq\RK{\calA}$.\rq\rq\  Moreover, the constituent of the $\ell^0$ minimization model used to recover the signal are also the constituents of \cref{model1}. Again, the main novelty is to extend this result to the identifiability of the factors of deep matrix products.

In order to establish a necessary condition for identifiability, first note that if we extend
$P(\Mod{L}) - P(\Mod{L'})$ to be scale invariant, this will not
affect whether or not it intersects $\ker(\calA)$ outside of the
origin. We immediately conclude that in the complex setting
where $\Mod{L},\Mod{L'}$ are both Zariski closed, that
$\Mod{}$ is non-identifiable whenever $\rk(\calA) < \DIM {P(\Mod{L}) - P(\Mod{L'})}$.
This indicates that we should always expect non-identifiability  whenever $\rk(\calA) < \DIM {P(\Mod{L}) - P(\Mod{L'})}$ but is not adequate to prove it because real algebraic varieties
need not satisfy (*). However it is true for real linear spaces, so
we immediately conclude  the following weak result:

\begin{thm}{\bf Necessary condition for Identifiability}\label{necess-ident-thm}

Let $C(P(\Mod{L}) - P(\Mod{L'}))$ be the set of all points on all lines through
 the origin intersecting $P(\Mod{L}) - P(\Mod{L'})$, and let
 $q$ be the maximal dimension of a linear space on $C(P(\Mod{L}) - P(\Mod{L'}))$.
Then if $q>\rk(\calA)$, $\Mod{}$ is not identifiable.
In particular when the $P(\Mod{L})$'s contain linear space and if we let $S'$ the be the largest dimension of these vector space, if $2 S' >\rk(\calA)$, then $\Mod{}$ is not identifiable.
\end{thm}

Let us illustrate the theorems by considering a deep feed-forward ReLU network and consider the structured linear network obtained by fixing the action of ReLU (as is done in \cref{struct-linear-network-sec}). The matrix $X$ contains the outputs, the operator $M_K$ multiplies the matrix containing the inputs by the weights between the first and the second layer. For every input/output pair, the action of ReLU is different and removes paths from the input entries to the output entries. We assume however that every entry of every output is reached by at least one path in the network starting at a non-zero entry of the input. In that case, it is not difficult to see that $\calA$ is a surjection and therefore $\rk(\calA) = mn$, where $m$ is the size of the output and $n$ is the number of learning samples.

The condition in \cref{suffic-ident-thm} becomes
\[mn \geq \DIM{\Sigma_2} = 2K(S-1) + 2,
\]
and $KS$ is typically the number of parameters of the network. The intuition behind \cref{suffic-ident-thm} is that, if the action of ReLU is sufficiently random and if the above inequality holds, we can expect the network to be identifiable with high probability\footnote{This statement gives the intuition behind \cref{suffic-ident-thm} but it should be made precise, as emphasized in the perspectives of this paper.}. This situation corresponds to an under-parameterized case (favorable for identifiability).

The condition in \cref{necess-ident-thm} is
\[2S > mn.
\]
When this inequality holds the network is not identifiable. It corresponds to an over-parameterized configuration.  In the intermediate situation, when $2S\leq mn < 2K(S-1) + 2$, and when the action of the activation function does not introduce sufficiently randomness, the theorems are inconclusive.

Notice that such networks can also be analyzed using \cref{ident-thm}. It is indeed not difficult to see that if there exists two paths that: 1/ start from the same entry of the input layer; 2/ end at the same entry of the output layer; 3/ if both paths are present (despite the action of ReLU) for every input/output pair; then \cref{ident-hS} does not hold\footnote{Simply consider two rank one tensors, each tensor being a Dirac at the position corresponding to one of the two paths.} and $ \RR^{S\times K}$ is not identifiable. It is not clear at this point that the conditions 1, 2, 3 are met by all non-identifiable structured linear feed-forward networks. However, removing paths from the network (as is done by ReLU and Dropout) is a way to avoid conditions 1, 2, 3 to be met.

%
\section{Stability guarantee}\label{stable-sec}

In this section, we consider errors of different natures. 
We assume that there exists $\overline{L}$ and $L^*\in\NN$, $\overline{\hbf}\in\Mod{\overline{L}}$ and $\hbf^*\in\Mod{L^*}$, such that
\begin{equation}\label{delta}
\|M_1(\overline{\hbf}_1) \cdots M_K(\overline{\hbf}_K)- X\| \leq \delta,
\end{equation}
and
\begin{equation}\label{inexact_model}
\|M_1(\hbf_1^*)\cdots M_K(\hbf_K^*) -X \| \leq \eta,
\end{equation}
for $\delta$ and $\eta$ typically small.

Again, this corresponds to existing unknown parameters $\overline \hbf$ that we estimate from a noisy observation $X$, using an inaccurate solution $\hbf^*$ of \cref{model} (as in \cite{bourrier2014fundamental} where the case $K=1$ is studied). Otherwise, $\overline \hbf$ and $\hbf^*$ shall be interpreted as different learned parameters; $\delta$ and $\eta$ are the corresponding risks.




Notice that the above hypothesis does not even require \cref{model} to have a solution. Algorithms which do not come with a guarantee sometimes manage to reach small $\delta$ and $\eta$ values. In those cases, the analysis we conduct in this section permits to get the stability guarantee, despite the lack of a  guarantee of the algorithm. Finally, the hypotheses \cref{delta}~ and \cref{inexact_model} enable one to obtain guarantees for algorithms that, instead of minimizing  \cref{model}, minimize an objective function which approximates the one in \cref{model}. This is particularly relevant for machine learning applications when \cref{model} can be an empirical risk that needs to be regularized or is not truely minimized (for instance, when using {\em dropout} \cite{srivastava2014dropout}).

A necessary and sufficient condition for the identifiability of $\Mod{}$ is stated in \cref{ident-thm}. The condition is on the way $\KER{\calA}$ and $P(\Mod{L}) - P(\Mod{L'})$ intersect. In order to get a stability guarantee, we need a stronger condition on the geometry of this intersection to hold for every $L$ and $L'\in\NN$. This condition is provided in the next definition.

\begin{defi}{\bf Deep-Null Space Property}\label{dnsp-def}

Let $\gamma>0$ and $\rho>0$. We say that $\KER{\calA}$ satisfies the {\em \NSPlong~(\NSP)} with respect to the collection of models $\Mod{}$ with constants $(\gamma,\rho)$ if for any $L$ and $L'\in\NN$, any $T\in P(\Mod{L}) - P(\Mod{L'})$ satisfying $\|\calA T\|\leq \rho$ and any $T'\in\KER{\calA}$, we have
\begin{equation}\label{dnsp}
\|T\| \leq \gamma \|T-T'\|.
\end{equation}
\end{defi}
The \NSP~implies that,  for $T\in P(\Mod{L}) - P(\Mod{L'})$ close to $\KER{\calA}$ in the sense that $\|\calA T\|\leq \rho$ we must have, by decomposing $T= T'+T''$, with $T'\in\KER{\calA}$ and $T''$ in its orthogonal complement
\[\|T\| \leq \gamma \|T-T'\| = \gamma \|T''\| \leq \frac{\gamma}{\sigma_{min}} \|\calA T''\| \leq \frac{\gamma}{\sigma_{min}} \rho,
\]
where $\sigma_{min}$ is the smallest non-zero singular value of $\calA$. In words, $\|T\|$ must be small. We can conclude that under the \NSP, $P(\Mod{L}) - P(\Mod{L'})$ and $\{T\in\TS \mid \|\calA T\|\leq \rho\}$ intersect at most in the neighborhood of $0$.

Additionally, \cref{dnsp} implies that in the neighborhood of $0$, $\KER{\calA}$ and $P(\Mod{L}) - P(\Mod{L'})$ are not tangential, i.e., their intersection is transverse.

If $\KER{\calA}$ satisfies the \NSP~with respect to the collection of models $\Mod{}$ with constants $(\gamma,\rho)$, then for all $T'\in\KER{\calA}$ and all $T\in P(\Mod{L}) - P(\Mod{L'})$ satisfying $\|\calA T\|\leq \rho$,
\[\|T'\| \leq \|T\| + \|T'-T\| \leq (\gamma +1) \|T'-T\|.
\]
Therefore,
\begin{equation}\label{usualNSP}
\forall T'\in\KER{\calA},\qquad \|T'\| \leq (\gamma +1) d_{loc}(T', P(\Mod{L}) - P(\Mod{L'})),
\end{equation}
where we have set for any $C\subset\TS$
\[d_{loc}(T',C) = \inf_{T\in C, \|\calA T\|\leq \rho}  \|T'-T\|.
\]
The converse is also true: if $\KER{\calA}$ satisfies \cref{usualNSP}, it satisfies the \NSP~with respect to the collection of models $\Mod{}$ with appropriate constants. In the context of the usual compressed sensing (i.e., when $K=1$, $\Mod{L}$ contains $L$-sparse signals, $\calA$ is a rectangular matrix with full row rank and $X$ is a vector), the localization appearing in $d_{loc}$ can be discarded since the inequality must hold when $T'$ is small and since in this case this localization has no effect. Therefore, in the compressed sensing context, \cref{usualNSP} (and therefore \NSP) is the usual Null Space Property with respect to $L$-sparse vectors, as defined in \cite{cohen2009compressed}. However, \NSP~is generalized to take into account deep structured linear network. This motivates the name.

In the general case, the \NSP~can be understood as a local version of the generalized-NSP for $\calA$ relative  to $P(\cup_{L\in\NN}\Mod{L}) - P(\cup_{L\in\NN}\Mod{L})$, as defined in \cite{bourrier2014fundamental}. Our interest in locality (as imposed by the constraint $\|\calA T\|\leq \rho$) is motivated by the fact that we want to use the \NSP~when the signal to noise ratio is controlled (i.e., the hypotheses of \cref{rk1-Ident-thm} are satisfied). The condition for the stability property therefore includes such hypotheses.

We have not adapted the robust-NSP defined in \cite{bourrier2014fundamental}. The benefit in not using this definition is to obtain a simpler definition for \NSP. In particular \cref{dnsp} does not involve the geometry of $\calA$ in the orthogonal complement of $\KER{\calA}$. Looking in detail at the benefit of this adaptation is of great interest.

Finally, we trivially have the following facts:
\begin{itemize}
\item If $\KER{\calA} = \{0\}$, then  $\KER{\calA}$ satisfies the \NSP~with respect to the model $\hS$ with constant $(1,+\infty)$.
\item For any $\gamma'\geq \gamma$: If $\KER{\calA}$ satisfies the \NSP~with respect to the collection of models $\Mod{}$ with constants $(\gamma,\rho)$, then $\KER{\calA}$ satisfies the \NSP~with respect to the collection of models $\Mod{}$ with constants $(\gamma',\rho)$.
\item For any $\Mod{'}\subset \Mod{}$: If $\KER{\calA}$ satisfies the \NSP~with respect to the collection of models $\Mod{}$ with constant $(\gamma,\rho)$, then $\KER{\calA}$ satisfies the \NSP~with respect to the collection of models $\Mod{'}$ with constants $(\gamma,\rho)$. In particular, if $\KER{\calA}$ satisfies the \NSP~with respect to the model $\hS$ with constant $(\gamma,\rho)$, it satisfies the \NSP~with respect to any collection of models, with constants $(\gamma,\rho)$.
\end{itemize}

\begin{thm}{\bf Sufficient condition for the stability property}\label{suf-stble-recovery}

Assume  $\KER{\calA}$ satisfies the \NSP~with respect to the collection of models $\Mod{}$ and with the constants $(\gamma,\rho)$. For any $\hbf^*$ as in \cref{inexact_model} with $\eta$ and $\delta$ (see \cref{inexact_model} and \cref{delta}) such that $\delta + \eta \leq \rho$, we have
\[ \|P(\hbf^*) - P(\overline\hbf)\|\leq \frac{\gamma}{ \sigma_{min}} ~(\delta+\eta),
\]
where  $\sigma_{min}$ is the smallest non-zero singular value of $\calA$. Moreover, if
$\overline{\hbf} \in\hS_*$ and $\frac{\gamma}{ \sigma_{min}} ~(\delta+\eta) \leq \frac{1}{2}  ~ \max\left(\|P(\overline\hbf)\|_\infty,\|P(\hbf^*)\|_\infty\right)$ then
\begin{equation}\label{tkl,stb}
d_p ( \CL{\hbf^*} , \CL{\overline\hbf} ) \leq \frac{7 (KS)^{\frac{1}{p}} \gamma}{\sigma_{min}}     \min\left( \|P(\overline\hbf)\|_\infty^{\frac{1}{K}-1}, \|P(\hbf^*)\|_\infty^{\frac{1}{K}-1} \right) (\delta+\eta).
\end{equation}
\end{thm}

The first part of the proof is very similar to standard proofs in the Compressed Sensing and stable recovery literature. The second part simply uses \cref{rk1-Ident-thm}. 
\ifthenelse{\boolean{long}}
{The theorem is proved in \cref{suf-stble-recovery-proof}.}
{The theorem is proved in the supplementary material \cref{suf-stble-recovery-proof} and in \cite{MalgouyresLandsbergLong}.
}


\cref{suf-stble-recovery} provides a sufficient condition to obtain stability. The only significant hypothesis made on the deep structured linear network is that $\KER{\calA}$ satisfies the \NSP~ with respect to the collection of models $\Mod{}$. One might ask whether this hypothesis is sharp or not. The next theorem shows that the answer to this question is positive.

\begin{thm}{\bf Necessary condition for the stability property}\label{nec-stble-recovery}

Assume the stability property holds: There exists  $C$ and $\delta>0$ such that for any $\overline{L}\in\NN$, $\overline{\hbf}\in\Mod{\overline{L}}$, any $X=\calA P(\overline{\hbf}) + e$, with $\|e\| \leq \delta$, any $L^*\in\NN$ and any $\hbf^*\in\Mod{L^*}$
such that
\[\|\calA P(\hbf^*) -X \|^2 \leq \| e \|
\]
we have
\[d_2 ( \CL{\hbf^*} , \CL{\overline\hbf} ) \leq  C ~ \min\left( \|P(\overline\hbf)\|_\infty^{\frac{1}{K}-1}, \|P(\hbf^*)\|_\infty^{\frac{1}{K}-1} \right)  \|e\|.
\]

Then, $\KER{\calA}$ satisfies the \NSP~ with respect to the collection of models $\Mod{}$ with constants $$(\gamma,\rho) =  (C S^{\frac{K-1}{2}} \sqrt{K}~\sigma_{max} , \delta)$$ 
where $\sigma_{max}$ is the spectral radius of $\calA$.
\end{thm}

The first part of the proof is inspired  by and similar to the proof of the analogous converse statement in \cite{cohen2009compressed}. The second part simply uses \cref{PLip-thm}. 
\ifthenelse{\boolean{long}}
{The theorem is proved in \cref{nec-stble-recovery-proof}.}
{The theorem is proved in the supplementary material \cref{nec-stble-recovery-proof} and in \cite{MalgouyresLandsbergLong}.
}

The sharpness of the known results when $K=2$ is usually argued by comparing the number of samples necessary for the recovery and the information theoretic limit of the problem. As far as the authors know, the above theorem is therefore new even when $K=2$.

As is usually the case with Null Space Property or Restricted Isometry Property, it will often be difficult or impossible to establish that a particular operator ${\calA}$ satisfies the \NSP~ with respect to the collection of models $\Mod{}$. To find favorable cases, we need to consider random operators $\calA$ such that the distribution of $\calA$ enables one to establish that the \NSP~ holds with high probability, when in the right configurations (see the bibliography in \cref{bibli-ref} whose references contain many examples of such arguments). The most common distribution for the analog of $\calA$ include operators/matrices whose coefficients are Gaussian or Bernoulli. Collections of models inducing sparsity, non-negativity, low-rank constraints {\it etc}. are the most studied. In this regard, the fact that there exists a low complexity test guaranteeing that the networks considered in \cref{conv-tree} can be stably recovered is an exception.


\section{Application to convolutional linear network}\label{conv-tree}

\begin{figure}
\centering{
\begin{picture}(15,10)
\put(1,10){\shortstack{edges of depth}}
\put(6,10){\shortstack{$3$}}
\multiput(6.2,9.7)(0,-0.1){90}{\line(0,-1){0.03}}

\put(10,10){\shortstack{$2$}}
\multiput(10.2,9.7)(0,-0.1){90}{\line(0,-1){0.03}}

\put(14,10){\shortstack{$1$}}
\multiput(14.2,9.7)(0,-0.1){90}{\line(0,-1){0.03}}

\put(0,4.8){\shortstack{leaves}}
\multiput(2,5)(0.1,0.3){10}{\line(1,3){0.07}}
\multiput(2,5)(0.1,0.1){10}{\line(1,1){0.07}}
\multiput(2,5)(0.1,-0.1){10}{\line(1,-1){0.07}}
\multiput(2,5)(0.1,-0.3){10}{\line(1,-3){0.07}}
\put(3,1.8){\shortstack{$\RP$}}
\put(4,2){\line(1,0.25){4}}\put(4,2){\vector(1,0.25){1.8}}
\put(3,3.8){\shortstack{$\RP$}}
\put(4,4){\line(1,-0.25){4}}\put(4,4){\vector(1,-0.25){1.8}}
\put(3,5.8){\shortstack{$\RP$}}
\put(4,6){\line(1,0.25){4}} \put(4,6){\vector(1,0.25){1.8}}
\put(3,7.8){\shortstack{$\RP$}}
\put(4,8){\line(1,-0.25){4}} \put(4,8){\vector(1,-0.25){1.8}}

\put(8,3){\line(1,0){4}}\put(8,3){\vector(1,0){1.8}}
\put(8,3){\line(1,1){4}}\put(8,3){\vector(1,1){1.8}}
\put(8,7){\line(1,0){4}}\put(8,7){\vector(1,0){1.8}}
\put(8,7){\line(1,-1){4}}\put(8,7){\vector(1,-1){1.8}}
\put(12,3){\line(1,0.5){4}}\put(12,3){\vector(1,0.5){1.8}}

\put(12,7){\line(1,-0.5){4}}\put(12,7){\vector(1,-0.5){1.8}}
\put(16.5,4.8){\shortstack{$\RP$}}
\put(18,4.8){\shortstack{root $r$}}

\end{picture}
}
\caption{\label{network}Example of the considered convolutional linear network. To every edge is attached a convolution kernel. The network does not involve non-linearities or sampling.}
\end{figure}
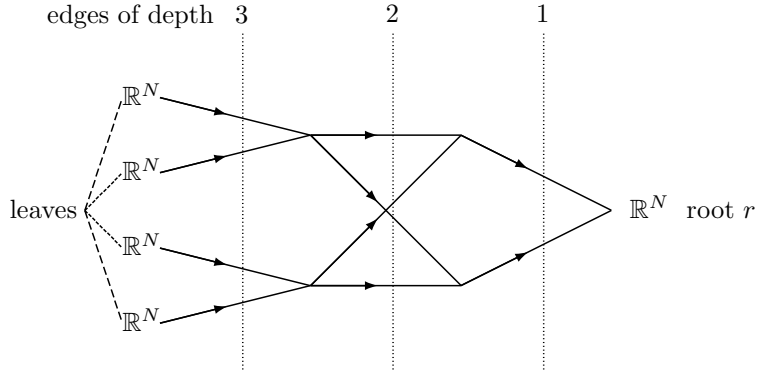

We consider a convolutional linear network as depicted in \cref{network}. The network typically aims at performing a linear analysis or synthesis of a signal living in $\RP$. The considered convolutional linear network is defined from a rooted directed acyclic graph $\tree(\edges,\nodes)$ composed of nodes $\nodes$ and edges $\edges$. Each edge connects two nodes. The root of the graph is denoted by $r$ and the set containing all its leaves is denoted by $\leaves$. We denote by $\PP$ the set of all paths connecting the leaves and the root. We assume, without loss of generality, that the length of any path between any leaf and the root is independent of the considered leaf and equal to some constant $K\geq 0$. We also assume that, for any edge $e\in\edges$, the number of edges separating $e$ and the root is the same for all paths between $e$ and $r$. It is called the depth of $e$. We also say that $e$ belongs to the layer $k$. For any $k\in\NNN{K}$, we denote the set containing all the edges of depth $k$, by $\edges(k)$. 

Moreover, to any edge $e$ is attached a convolution kernel of support $\calS_e\subset \NNN{N}$. We assume (without loss of generality) that $\sum_{e\in\edges(k)} |\calS_e|$ is independent of $k$ ($|\calS_e|$ denotes the cardinality of $\calS_e$). We take
\[S= \sum_{e\in\edges(1)} |\calS_e|.
\]
For any edge $e$, we consider the mapping $\calT_e: \RR^S \longrightarrow \RR^N$ that maps any $h\in\RR^S$ into the convolution kernel $h_e$, attached to the edge $e$, whose support is $\calS_e$. It simply writes at the right location (i.e., those in $\calS_e$) the entries of $h$ defining the kernel on the edge $e$. 

At each layer $k$, the convolutional linear network computes, for all $e\in\edges(k)$,  the convolution between the signal at the origin of $e$; then, it attaches to any ending node the sum of all the convolutions arriving at that node. Examples of such convolutional linear networks includes wavelets, wavelet packets \cite{Mallatbook} or the fast transforms optimized in \cite{FTL_IJCV,FTO}. It is clear that the operation performed at any layer depends linearly on the parameters $h\in\RR^S$ and that its results serve as inputs for the next layer. The convolutional linear network therefore depends on parameters $\hbf\in\hS$ and takes the form
\[X = M_1(\hbf_1)\cdots M_K(\hbf_K),
\]
where the operators $M_k$ satisfy \cref{mkdef}. 

This section aims at identifying conditions such that any unknown parameters $\overline{\hbf}\in\hS$ can be identified or stably recovered from $X = M_1(\overline\hbf_1)\cdots M_K(\overline\hbf_K)$ (possibly corrupted by an error).

In order to do so, we introduce some notation. We apply the convolutional linear network to an input $x\in\codeset$, where $x$ is the concatenation of the signals $x^f\in\RR^N$ for $f\in\leaves$. Therefore, $X$ is the (horizontal) concatenation of $|\leaves|$ matrices $X^f\in\RR^{N\times N}$ such that
\begin{equation}\label{revnoiritb}
X x = \sum_{f\in\leaves} X^fx^f\qquad\mbox{, for all }x\in\codeset.
\end{equation}
Consider the convolutional linear network defined by $\hbf\in\hS$ as well as $f\in\leaves$ and $n\in\NNN{N}$. The column of $X$ corresponding to the entry $n$ in the leaf $f$ is the translation by $n$ of
\begin{equation}\label{multiconv}
\sum_{p\in\PP(f)} \multiconv{p}{\hbf}
\end{equation}
where $\PP(f)$ contains all the paths of $\PP$ starting from the leaf $f$ and
\[\multiconv{p}{\hbf} = \calT_{e^1}(\hbf_1)*\cdots*\calT_{e^K}(\hbf_K) \qquad\mbox{, with }p=(e^1,\cdots,e^K),
\]
is the composition of convolutions along the path $p$.

For any $k\in\NNN{K}$, define the mapping $\ebf_k: \NNN{S} \longrightarrow \edges(k)$ which provides for any $i\in\NNN{S}$ the unique edge of $\edges(k)$ such that the $i^{\mbox{th}}$ entry of $h\in\RR^S$ contributes to $\calT_{\ebf_k(i)}(h)$. For any $\ibf\in\NNN{S}^K$, let $\pbf_\ibf=(\ebf_1(\ibf_1), \cdots, \ebf_K(\ibf_K))$ and
\[\Ibf = \left\{\ibf \in\iS | \pbf_\ibf \in \PP \right\}.
\] 
The latter contains all the indices corresponding to a valid path in the network. For any set of parameters $\hbf\in\hS$ and any path $\pbf\in\PP$, we also let $\hbf^\pbf$ denote  the restriction of $\hbf$ to its indices contributing to the kernels on the path $\pbf$. We let $\one\in\RR^S$ denote a vector of size $S$ with all its entries equal to $1$. For any edge $e$, $\one^e\in\RR^S$ consists of zeros except for the entries corresponding to the edge $e$ which are equal to $1$. For any $\pbf=(e^1,\cdots,e^K)\in\PP$, the support of $M_1(\one^{e_1})\cdots M_K(\one^{e_K})$ is denoted by $\SUPP{\pbf}$.

Finally, by \cref{AP=M-cor} there exists a unique mapping 
\[\calA:\TS\longrightarrow \RR^{N\times N|\leaves|}\] 
such that 
\[\calA P(\hbf) = M_1(\hbf_1)\cdots M_K(\hbf_K)\qquad \mbox{, for all } \hbf\in\hS,
\]
where $P$ is the Segre embedding defined in \cref{defP}. 

\begin{defi}
We say the topology of the network is {\em sufficiently scattered} if and only if  all the entries of $M_1(\one)\cdots M_K(\one)$ belong to $\{0,1\}$.
\end{defi}

The following statements will show that having a sufficiently scattered topology is a necessary and sufficient condition for the stability of the optimal parameters. Before going into this, we illustrate the scattering property with a simple example.

Consider a simple composition of two convolutions ($K=2$ and $|\PP| =1$). At first, we make an assumption on the supports $\calS_e$, imposing that  the supports of both kernels are in $\{1,2,3\}$. The topology is obviously not sufficiently scattered. Indeed, some of the entries of the convolution kernel corresponding to the matrix $M_1(\one)M_2(\one)$ are equal to $2$.

Now consider an assumption on the network topology imposing $\{1,2,3\}$ for the support of the first kernel and $\{1,10\}$ for the second. When we observe the convolution of two kernels having such supports, we see two replicas of the first kernel; the amplitudes of the replicas depend on the second kernel; and both kernels are identifiable. In this last example, the network topology is sufficiently scattered.

The scattering condition can easily be computed using \cref{algo_check}. Indeed, when applying the network to a Dirac in the leaf $f$, using \cref{revnoiritb}, we obtain the convolution kernel of $X_f$. We can then easily test if $X_f$ only contains $0$'s and $1$'s. The numerical complexity of \cref{algo_check} is essentially the cost for applying $|\leaves|$ times the network. It is usually low. Notice that a network is sufficiently scattered if and only if, for all leaves $f\in\leaves$,  the sub-networks originating at $f$ are sufficiently scattered. The scattering of these subnetworks are independent. The fact the the convolution kernels, for the different leaves, overlap does not affect the scattering property.

\begin{algorithm}                   
\begin{algorithmic}              
\REQUIRE The network topology
\ENSURE Boolean output = \lq true\rq, if the topology is sufficiently scattered;  \lq false\rq, otherwise.
\STATE 
\STATE output = true
\FOR{ $f\in\leaves$}
\STATE Build $x$: a Dirac positioned at the leaf $f$
\STATE Apply the network to $x$ in order to compute $y=M_1(\one)\cdots M_K(\one) x$
\STATE If some of the entries of $y$ are outside $\{0,1\}$ then set output = false
\ENDFOR
\end{algorithmic}
\caption{Algorithm testing if the topology of the convolutional network leads to the stability guarantee. \label{algo_check}}                   
\end{algorithm}

Finally, beside the known examples in blind-deconvolution (i.e., when $K=2$ and $|\PP|=1$) \cite{ahmed2014blind,bahmani2015lifting}, there are (truly deep) convolutional linear networks satisfying the condition of the first statement of  \cref{nec-ident-network}. For instance, the convolutional linear network corresponding to the un-decimated Haar (wavelet)\footnote{Un-decimated means computed with the \lq\lq Algorithme \`a trous\rq\rq, \cite{Mallatbook}, Section 5.5.2 and 6.3.2. The Haar wavelet is described in \cite{Mallatbook}, Section 7.2.2, p. 247 and Example 7.7, p. 235.} transform is a tree and for any of its leaves $f\in\leaves$, $|\PP(f)|=1$. Moreover, the support of the kernel living on the edge $e$, of depth $k$, on this path is $\{0, 2^k\}$. It is therefore not difficult to check that the first condition of \cref{nec-ident-network} holds.

\begin{prop}{\bf Necessary condition of identifiability of convolutional linear network}\label{nec-ident-network}

\begin{itemize}
\item Either the topology of the network is sufficiently scattered and then
\begin{enumerate}
\item \label{premier_item} for any distinct $\pbf$ and $\pbf'\in\PP$, $\SUPP{\pbf} \cap \SUPP{\pbf'}  = \emptyset$.
\item $\KER{\calA} = \{T\in\TS | \forall \ibf\in\Ibf, T_\ibf = 0\}.$
\end{enumerate}
\item or the topology of the network is not sufficiently scattered and then $\hS$ is not identifiable.
\end{itemize}

\end{prop}

\ifthenelse{\boolean{long}}
{The Proposition is proved in \cref{nec-ident-network-proof}.}
{The Proposition is proved in the supplementary material and in \cite{MalgouyresLandsbergLong}.
}

\begin{prop}\label{network-cor}
If $|\PP|=1$ and the topology of the network is sufficiently scattered, then $\KER{\calA} =\{0\}$ and $\KER{\calA}$ satisfies the \NSP~with respect to any model collection $\Mod{}$ with constant $(\gamma,\rho)=(1,+\infty)$. Moreover, we have $\sigma_{min} = \sqrt{N}$. 
\end{prop}

\ifthenelse{\boolean{long}}
{The Proposition is proved in \cref{network-cor-proof}.}
{The Proposition is proved in the supplementary material and in \cite{MalgouyresLandsbergLong}.
}

In what follows, we establish stability results for a convolutional linear network estimator. In order to do so, we consider a convolutional linear network of known structure $\tree(\edges,\nodes)$ and $(\calS_e)_{e\in\edges}$. We consider parameters $\overline \hbf \in\hS$ and $\hbf^* \in\hS$ such that
\begin{equation}\label{eq-hbar}
\| M_1(\overline\hbf_1)\cdots M_K(\overline\hbf_K) -X \| \leq \delta,
\end{equation} 
and
\begin{equation}\label{eq-hstar}
\| M_1(\hbf^*_1)\cdots M_K(\hbf^*_K) - X \|\leq\eta.
\end{equation} 
We say that two networks sharing the same structure and defined by $\hbf$ and $\gbf\in\hS$ are equivalent if and only if
\[\forall \pbf\in\PP, \exists (\lambda_e)_{e\in\pbf} \in\RR^\pbf\mbox{, such that } \prod_{e\in\pbf} \lambda_e = 1 \mbox{ and } \forall e\in\pbf, \calT_e(\gbf) = \lambda_e \calT_e(\hbf). 
\]
The equivalence class of $\hbf\in\hS$ is denoted by $\class{\hbf}$. For any $p\in[1,+\infty]$, we define
\[\cald_p (\class{\hbf}, \class{\gbf}) = \left(\sum_{\pbf\in\PP} d_p(\CL{\hbf^\pbf},\CL{\gbf^\pbf})^p  \right)^{\frac{1}{p}},
\]
where recall that $\hbf^\pbf$ (resp $\gbf^\pbf$) denotes the restriction of $\hbf$ (resp. $\gbf$) to the path $\pbf$ and $d_p$ is defined in \cref{def_d}. Since $d_p$ is a metric, it follows that $\cald_p$ is a metric between network classes.

We summarize the results concerning convolutional networks in the following theorem.

\begin{thm}{\bf Necessary and sufficient condition of stable recovery of convolutional linear network}\label{stabl-rec-network-thm}

If \cref{algo_check} returns \lq false\rq, the network topology is not sufficiently scattered and the network is not identifiable.

If \cref{algo_check} returns \lq true\rq, if $\overline\hbf$ and $\hbf^*$ satisfy \cref{eq-hbar} and \cref{eq-hstar} and 
\begin{itemize}
\item if all the edges support a significant convolution kernel: there exists $\varepsilon >0$ such that for all $e\in\edges$, $\|\calT_e(\overline\hbf)\|_\infty\geq \varepsilon$, 
\item if the \lq\lq signal to noise ratio\rq\rq is sufficient: $\delta+\eta \leq \frac{\sqrt{N} \varepsilon^K}{2}$, 
\end{itemize}
then the network defined by $\hbf^*$ and $\overline\hbf$ are close to each other
\[\cald_p (\class{\hbf^*}, \class{\overline\hbf}) \leq 7 (KS')^{\frac{1}{p}} \varepsilon^{1-K}~ \frac{\delta+\eta}{\sqrt{N}},
\]
where $S' = \max_{e\in\edges} |\calS_e|$ is the size the largest convolution kernel.
\end{thm}

\ifthenelse{\boolean{long}}
{The theorem is proved in \cref{stabl-rec-network-thm-proof}.}
{The theorem is proved in the supplementary material \cref{stabl-rec-network-thm-proof} and in \cite{MalgouyresLandsbergLong}.
}

\section{Conclusion and perspectives}

In this paper, we have established necessary and sufficient conditions for the identifiability and stable recovery of deep structured linear networks. They rely on the lifting of the problem in a tensor space. The technique is called {\em tensorial lifting}. The main results are proved using compressed sensing technics  and properties of the Segre embedding (the embedding that maps the parameters in the tensor space). The general results are then specialized to establish necessary and sufficient conditions for the stable recovery of a convolutional linear network of any depth $K\geq 1$.

Among the most salient perspectives, we mention the possibility to study deep feed-forward ReLU networks. For such a network, the action of ReLU is different for every sample; this leads to a different operator ${\mathcal A}$ for every sample; and all the different ${\mathcal A}$'s sense (linearly) the same rank one tensor. We can concatenate these operators to form a unique sensing operator. For instance, when modeling the action of ReLU as a Bernouilli variable applied to every path of the network, we expect to obtain sample complexity bounds (for instance) under the favorable hypothesis that an oracle has given us the action of ReLU.

A natural perspective of this work is also to study compressed networks (see \cite{arora2018stronger}), when the compression preserves the expressivity of the network.

Finally, the model considered in this paper approximately solves polynomial equations: $\calA P(\hbf) \sim X$. The structure of the polynomials is induced by the operators $M_k$ (i.e. the network topology) is very particular and restrictive. For instance, we only consider homogene polynomials in $P(\hbf)$. Extending this work to larger families of polynomials as well as limits of polynomials seem natural.


\ifthenelse{\boolean{long}}
{\section{Appendices}}
{
}

\ifthenelse{\boolean{long}}
{\PROOFUN}
{
}

\ifthenelse{\boolean{long}}
{\PROOFDEUX}
{
}

\ifthenelse{\boolean{long}}
{\PROOFTROIS}
{
} 


\ifthenelse{\boolean{long}}
{\PROOFQUATRE}
{
} 


\ifthenelse{\boolean{long}}
{\PROOFCINQ}
{
}


\ifthenelse{\boolean{long}}
{\PROOFSIX}
{
}

\ifthenelse{\boolean{long}}
{\PROOFSEPT}
{
}


\ifthenelse{\boolean{long}}
{\PROOFHUIT}
{
}


\ifthenelse{\boolean{long}}
{\PROOFNINE}
{
}


\ifthenelse{\boolean{long}}
{\PROOFDIX}
{
}


\ifthenelse{\boolean{long}}
{\PROOFONZE}
{
}


\ifthenelse{\boolean{long}}
{\PROOFDOUZE}
{
}


\ifthenelse{\boolean{long}}
{\PROOFTREISE}
{
}


\ifthenelse{\boolean{long}}
{\PROOFQUATORZE}
{
}


\bibliographystyle{plain}
\bibliography{ref}

\end{document}


\maketitle
\section{Proofs}
\PROOFUN

\PROOFDEUX

\PROOFTROIS

\PROOFQUATRE

\PROOFCINQ

\PROOFSIX   

\PROOFSEPT

\PROOFHUIT

\PROOFNINE

\PROOFDIX

\PROOFONZE

\PROOFDOUZE

\PROOFTREISE

\PROOFQUATORZE

